\documentclass{amsart}

\usepackage{amsmath,amssymb}

\newcommand{\R}{\mathbb{R}}
\newcommand{\C}{\mathbb{C}}
\newcommand{\Z}{\mathbb{Z}}
\newcommand{\N}{\mathbb{N}}
\renewcommand{\H}{\mathbb{H}}

\newcommand{\pol}{\textnormal{pol}}

\newcommand{\g}{{\mathfrak{g}}}    
\newcommand{\h}{{\mathfrak{h}}}    
\newcommand{\n}{{\mathfrak{n}}}

\renewcommand{\L}{{\mathcal{L}}}

\newcommand{\sA}{{\mathcal{A}}}
\newcommand{\F}{{\mathcal{F}}}

\newcommand{\ord}{\textnormal{ord}}
\newcommand{\Ad}{\textnormal{Ad}}

\newtheorem{Theorem}{Theorem}[section]
\newtheorem{Lemma}[Theorem]{Lemma}
\newtheorem{Def}[Theorem]{Definition}
\newtheorem{Cor}[Theorem]{Corollary}
\newtheorem{Prop}[Theorem]{Proposition}

{
{\theoremstyle{remark}\newtheorem{Remark}[Theorem]{Remark}}

\begin{document}
\setlength{\textheight}{209mm}

%
%
\title[Character formula]
{A character formula for representations of loop groups based on
  non-simply connected Lie groups}
\date{\today}
\author[R.~Wendt]{Robert Wendt}


\maketitle
\vspace{-.8cm}
\begin{center}
University of Toronto\\
Department of Mathematics\\
100 St. George Street\\
Toronto, Ontario, M5S 3G3\\
Canada\\

\medskip

E-mail: rwendt@math.toronto.edu

\medskip

\end{center}

\section{Introduction}
In this paper we compute characters of certain irreducible representations of
loop groups based on non simply connected Lie groups. Apart from a
general representation theoretic interest, our motivation to study these
characters comes from the fact that they appear naturally in the theory of
moduli spaces of semistable principal bundles over elliptic curves.

The characters of highest weight representations of loop
groups based on 
simply connected Lie groups are well understood due to the Kac-Weyl
character formula: Let $G$ be a simply connected complex Lie group whose 
Lie algebra is simple. We
denote by $L(G)$ the group of holomorphic maps from $\C^*$ to $G$. This
group possesses a universal central extension $\widehat L(G)$ which, 
viewed as a manifold, 
is a non-trivial $\C^*$-bundle over $L(G)$. The natural multiplicative
action of $\C^*$ on $L(G)$ lifts uniquely to an action of $\C^*$ on 
$\widehat L(G)$ by group automorphisms.
Let $V$ be an
irreducible highest weight representation of $\widehat L(G)$. Such a
representation extends to a representation of the semidirect product
$\widehat L(G)\rtimes \C^*$.
One can show that for any $q$ with $|q|<1$, the element 
$(g,q)\in \widehat L(G)\rtimes \C^*$ viewed as an operator on $V$
extends to a trace class operator
on the Hilbert space completion of $V$. So one can define the
character $\chi_V$ of $V$ at a point $(g,q)$ with $|q|<1$ 
as the trace of the
operator $(g,q)$. This defines a  holomorphic and conjugacy-invariant 
function on the space
$\widehat L(G)\times D^*$, where $D^*$ denotes the punctured unit disk
in $\C$. The Kac-Weyl character formula gives an explicit formula for
the character $\chi_V$ restricted to a certain family of tori in 
$\widehat L(G)\rtimes D^*$ in terms of theta functions. 
This is enough to describe the character
completely since almost every conjugacy class in $\widehat L(G)\times
D^*$ intersects this family of tori.

If the group $G$ is not simply connected, the loop group $L(G)$ 
consists of several connected components which are labeled by the
fundamental group of $G$. In this case, central extensions of
$L(G)$ have been constructed in \cite{T}. 
We will review this
construction in section \ref{prelim}.
Let $\widehat L(G)$ denote
such a central extension. The natural action of $\C^*$
on $L(G)$ does not lift to $\widehat L(G)$. Instead, a finite covering
$\widetilde\C^*$ of $\C^*$ acts on $\widehat L(G)$ covering the
natural $\C^*$-action on $L(G)$. So similar to the simply connected
case, we can consider the semidirect product $\widehat
L(G)\rtimes\widetilde\C^*$.
We are interested in 
representations of $\widehat L(G)\rtimes\widetilde\C^*$ which,
restricted to the connected component of 
$\widehat L(G)\rtimes\widetilde\C^*$ 
containing the identity, decompose into a
direct sum of irreducible highest weight representations. These
representations have been classified in \cite{T}.
Let $V$ be
such a representation. 
For  $\widetilde q\in\widetilde \C^*$, let $q$ denote the image
under the natural projection $\widetilde \C^*\to\C^*$.
As in the
simply connected case, one shows that any 
$(g,\widetilde q)\in \widehat L(G)\rtimes\widetilde\C^*$ with $\vert
q\vert<1$ extends to a trace class operator on the Hilbert space
completion of $V$. Thus, one can define the character $\chi_V$ of the
representation $V$ exactly as in
the simply connected case. 

The main goal of this paper is to give an
explicit formula for the character $\chi_V$ restricted to the connected
components of $\widehat L(G)\rtimes\widetilde\C^*$
which do not contain the identity element. This gives a generalization
of the Kac-Weyl character formula. 
While the usual approach to the
Kac-Weyl formula is 
rather algebraic, we work in a completely geometric setting. 
In
particular, we identify the characters with sections in certain line
bundles over an Abelian variety. To do this, we have to realize the
set of semisimple conjugacy classes  in a connected component of 
$\widehat L(G)\rtimes\widetilde D^*$ as the
total space of a line bundle over a family of Abelian varieties over
$D^*$. 
Then we show that the characters have
to satisfy a certain differential equation.
In the simply connected case, the
differential equation has been derived in \cite{EK}. Our main step is
a generalization of this equation to the non-simply connected case.
Finally, we use the differential equation to obtain an
explicit formula for the character $\chi_V$. In the simply connected
case, this gives an easy proof of the Kac-Weyl character formula which
is similar to Weyl's original proof of his character formula for
compact Lie groups. 
In the non-simply connected case, we obtain a formula for the
characters which very much resembles the 
Kac-Weyl character formula (Theorem \ref{KW2}). 
The main difference is that the character
restricted to a connected component of 
$\widehat L(G)\rtimes\widetilde\C^*$ not containing the identity is
not governed by the root system $\widetilde\Delta$ of the the Lie algebra of
$\widehat L(G)\rtimes\widetilde\C^*$ but by a new root system
$\widetilde\Delta_{\sigma_c}$ which can be obtained from $\widetilde\Delta$ by
a ''folding'' process. 
It is interesting to note that the Lie algebra
corresponding to root system $\widetilde\Delta_{\sigma_c}$ 
can, in general, not be realized
as a subalgebra of the affine Lie algebra corresponding to
$\widetilde\Delta$. In this way, the situation resembles the
case of characters of irreducible representations of non-connected
compact Lie groups \cite{W}. Also, Fuchs et al. \cite{FRS}, \cite{FSS}
have obtained similar results calculating the characters 
of representations of Kac-Moody
algebras twisted by outer automorphisms. These so called ''twining
characters'' appear in a
conjecture concerning Verlinde formulas for non-simply connected Lie
groups \cite{FS}.

Our main motivation for the study of characters of irreducible
representations of 
loop groups based on non simply connected Lie groups comes from the
theory of moduli spaces of semistable $G$-bundles on elliptic
curves. For a given group $G$ and an elliptic curve $E_q=\C^*/q^\Z$
with $q\in D^*$,
the moduli space of semistable $G$-bundles over $E_q$ consists of
several connected components which are labeled by the elements of the
fundamental group of $G$. The knowledge of the characters of
$\widehat L(G)\rtimes\widetilde\C^*$ allows to construct an analogue
of a Steinberg cross section in each connected component of 
$\widehat L(G)\times\{\widetilde q\}$ for any $\widetilde q\in\C^*$ such that
$q\in D^*$ 
(see e.g. \cite{Br} for the construction of
a Steinberg cross section in loop
groups based on simply connected Lie groups and \cite{M} for the case
of non-connected semisimple algebraic groups). It 
turns out that there is a natural action of $\C^*$ on this cross
section, and that the space of orbits of this action is isomorphic to
the connected component of the moduli space of semistable $G$--bundles 
on $E_q$
which corresponds to a connected component of 
$\widehat L(G)\times\{\widetilde q\}$. On the other hand, by
construction, the cross section can be identified with an 
affine space $\C^r$ and the $\C^*$--action becomes linear in this
identification. 
So the approach outlined above gives a new proof 
of a result of
Friedman and Morgan \cite{FM2} which states that each component
of the moduli space is isomorphic to  a weighted projective
space. These ideas will be published elsewhere.


The organization of this paper is as follows. In section \ref{prelim},
we review the construction of central
extensions $\widehat L(G)$ of loop groups $L(G)$ 
based on non simply connected Lie groups $G$ and
describe their representation theory.
In section \ref{conjugacy} we study certain conjugacy classes in
these groups.
In section \ref{proof} we identify the characters of the
representations of $\widehat L(G)$ introduced in section \ref{prelim}
with sections of a line bundle over a family of Abelian
varieties and deduce a differential equation for these
sections. Finally, in section \ref{char}, we put everything
together and deduce an explicit formula for the characters. In
particular, in \ref{sec:del} we describe how the ''folded'' root
system $\widetilde\Delta_{\sigma_c}$ appears for these characters. 
In the appendix we list the root systems
$\widetilde\Delta_{\sigma_c}$ together with some other data corresponding to
non simply connected Lie groups.



\section{Affine Lie groups and algebras}\label{prelim}


\subsection{Affine Lie algebras}\label{aff}
We begin by recalling some facts from the theory of affine Lie algebras.
Let $\g$ be a complex finite dimensional simple Lie algebra and 
let $\h\subset\g$ be a Cartan subalgebra. We denote the root 
system of $\g$ with respect to $\h$ by $\Delta$ and let $\h_\R\subset\h$ be 
the real vector space spanned by the co-roots of $\g$.  

The loop algebra $L(\g)$ of $\g$ is the Lie algebra of holomorphic 
maps from $\C^*$ to $\g$. The (untwisted) affine Lie algebra  
corresponding to $\g$ is a certain extension of $L(\g)$. Let us fix
some $k\in\C$ and 
consider the Lie algebra $\widetilde L^k(\g) = L\g\oplus\C C\oplus\C D$,
where the Lie bracket on $\widetilde L^k(\g)$ is given by
\begin{equation*} 
  [C,x(z)] = [C,D] = 0, \qquad [D,x(z)]=z\frac{d}{d z}x(z),
\end{equation*}
and
\begin{equation*}\label{etxeq}
  [x(z),y(z)]=[x,y](z)+\frac{k}{2\pi i}
  \int_{|z|=1}\langle \frac{d}{d z}x(z),y(z)\rangle dz\cdot C\,
\end{equation*}
%
%
Here $[x,y](z)$ denotes the pointwise commutator of $x$ and $y$, and 
$\langle.,.\rangle$ is the normalized invariant bilinear form on $\g$ (i.e.
the Killing form on $\g$ normalized in such a way that 
$\langle\alpha,\alpha\rangle=2$ for the long roots $\alpha$ of $\g$).
Note that the Lie algebras $\widetilde L^k(\g)$ are isomorphic for all
$k\neq0$. However, for different $k$ they define non-equivalent 
central extensions of $L(\g)\oplus\C D$. For $k=1$, we usually omit
the $k$ and denote the corresponding Lie algebra simply by
$\widetilde L(\g)$.

If $\g$ is simple, the subalgebra
$\widetilde L(\g)_{\pol}=\g\otimes\C[z,z^{-1}]\oplus\C C\oplus\C D
\subset\widetilde L(\g)$ of polynomial loops is an untwisted 
affine 
Lie algebra in the sense of  \cite{K}, and $\widetilde L(\g)$ 
can be viewed as 
a certain completion of it (see \cite{GW}). The Lie algebra 
$\widetilde L(\g)_{\pol}$ has a root space decomposition in the 
following sense: Set $\widetilde\h=\h\oplus\C C\oplus\C D$ and 
choose an element $\delta\in(\h_\R\oplus\R C\oplus\R D)^*$ dual 
to $D$. Then the root system $\widetilde\Delta$ of 
$\widetilde L(\g)_\pol$ is given by 
$$
  \widetilde\Delta=\{\alpha+n\delta~|~\alpha\in\Delta,~n\in\Z\}
   \cup\{n\delta~|~n\in\Z\setminus 0\}~
$$
and we can write
$$
  \widetilde L(\g) _\pol = 
   \widetilde\h\oplus\bigoplus_{\widetilde\alpha\in\widetilde\Delta}
     \widetilde L(\g)_{\widetilde\alpha}
$$
with $\widetilde L(\g)_{\widetilde\alpha}=\g_\alpha\otimes z^n$ if 
$\widetilde\alpha=\alpha+n\delta$, and 
$\g_{\widetilde\alpha}=\h\otimes z^n$ if $\alpha=n\delta$.


The set $\widetilde\Delta$ is an affine root system. 
Let $\Pi=\{\alpha_1,\ldots\alpha_l\}$ be a basis of $\Delta$ and let 
$\theta$ denote the highest root of $\Delta$ with respect to this basis. 
Then we can define 
$\widetilde\Pi=\{\alpha_0=\delta-\theta,\alpha_1,\ldots,\alpha_l\}$, 
which is a basis of $\widetilde\Delta$. The Dynkin diagram of 
$\widetilde\Delta$ is defined in the usual sense and it turns 
out that it is the extended Dynkin diagram corresponding to 
$\Delta$ (see \cite{K}). The affine root system $\widetilde\Delta$
decomposes into 
$$
  \widetilde\Delta =  \widetilde\Delta_+\cup\widetilde\Delta_-\,,
$$
where the set $\widetilde\Delta_+$ of positive roots is given by 
$$
  \widetilde\Delta_+ =
  \Delta_+\cup\{\alpha+n\delta~|~\alpha\in\Delta\cup\{0\},~n>0\}\,,
$$
and $\widetilde \Delta_-=-\widetilde \Delta_+$.

By definition, the set of real roots of $\widetilde L(\g)$ is the set
$$
  \widetilde\Delta^{re} = 
    \{\alpha+n\delta~|~\alpha\in\Delta\}\subset\widetilde\Delta\,,
$$
and the set of positive real roots is given by
$\widetilde\Delta_+^{re}=\widetilde\Delta^{re}\cap\widetilde\Delta_+$.

Sometimes we will need to consider twisted affine Lie algebras. 
If the finite dimensional Lie algebra $\g$ admits an outer automorphism 
$\sigma$ of finite order $\text{ord}(\sigma)=r$, 
one can define the twisted loop algebra
$$L(\g,\sigma)=\{X\in L(\g)~|~\sigma(X(z))=X(e^{\frac2r\pi i}z)\}~.$$ 
The corresponding affine Lie algebra $\widetilde L(\g,\sigma)$ is 
constructed in a similar manner as the untwisted algebra. 
It has a root space decomposition which is only slightly more 
complicated than in the untwisted case (see e.g. \cite{K}).


\subsection{Loop groups and affine Lie groups}
Let $G$ be a complex simply connected semisimple Lie 
group with Lie algebra $\g$ and suppose
that $\g$ is simple. The loop group $L(G)$ of $G$ is the group of 
holomorphic maps from  $\C^*$ to $G$ with pointwise 
multiplication. This is a Lie group with Lie algebra $L(\g)$.
Let $\widehat L(G)$
denote the universal central extension of $L(G)$. The central
extension $\widehat L(G)$ can be 
defined via the embedding of $L(G)$ into the ``differentiable loop''
group studied by Pressley and Segal \cite{PS}. 
Topologically, $\widehat L(G)$ is a non-trivial
holomorphic principal $\C^*$-bundle over $L(G)$.
In fact, there exists a central extension 
$\widehat L^k(G)$
of $L(G)$
for each $k\in\N$. The group $\widehat L^k(G)$ 
is called the level $k$ central extension of
$L(G)$. Its Lie algebra is given by
$\widehat L^k(\g) = \widetilde L^k(\g)/\C D$.
The universal central extension is just the level $1$ extension of
$L(G)$. 
The group $\C^*$ acts naturally on $L(G)$ by 
$(q\circ g)(z)=g(q^{-1}z)$ and we can consider the semidirect product
$L(G)\rtimes\C^*$.  
There is a  $\C^*$-action on $\widehat L^k(G)$ which covers the $\C^*$-action
on $LG$, and we  denote the semidirect product $\widehat L^k(G) \rtimes\C^*$ 
by $\widetilde L^k(G)$. Its Lie algebra is the affine Lie algebra 
$\widetilde L^k(\g)$ described in the last section

\medskip

Now assume  that $G$ is of the form $G=\widetilde G/Z$, where
$\widetilde G$ is simply connected and simple, and $Z\subset \widetilde G$ 
is a subgroup of the center of $\widetilde G$. 
Since the group 
$Z$ may be identified with the
fundamental group of $G$, the loop group $L(G)$ consists of 
$|Z|$
connected components.  In particular, the connected 
component of $L(G)$ containing the
identity element is isomorphic to $L(\widetilde G)/Z$.
We shall now indicate,
following Toledano Laredo \cite{T},  how
to construct certain central extensions of $L(G)$. We will first
consider the group $L_Z(\widetilde G)$ of holomorphic maps 
$g:\C\to \widetilde G$ such
that $g(t)g(t+1)^{-1}\in Z$. 
Identifying the variable $z$ with $e^{2\pi i t}$,
we see that the group $L(G)$ is isomorphic to $L_Z(\widetilde G)/Z$. 
Furthermore, the connected component of $L_Z(\widetilde G)$ containing the
identity element is isomorphic to $L(\widetilde G)$.

The goal is to
construct all central extensions of $L_Z(\widetilde G)$ 
and then see which of these
extensions are pullbacks of central extensions of $L(G)$. To this end, let 
$\widetilde T\subset\widetilde G$ and 
$T=\widetilde T/Z\subset G$ 
denote maximal tori of $\widetilde G$ and $G$, 
and let $\Lambda(\widetilde T)=Hom_{\text{alg grp}}(\C^*,\widetilde T)$ 
and $\Lambda(T)=Hom_{\text{alg grp}}(\C^*,T)$ denote
the respective co-character lattices. 
Then $\Lambda(T)/\Lambda(\widetilde T)\cong Z$.
The lattice $\Lambda(T)$ can be identified with a subgroup of
$L_Z(\widetilde G)$ by viewing it as a lattice in $\h_\R\subset\h$ and
identifying an element $\beta\in \Lambda(T)$ with the
''open loop'' $t\mapsto \exp(2\pi it\beta)$. We can define a subgroup 
$N\subset L(\widetilde G)\rtimes\Lambda(T)$ via 
$$
  N=\{(\lambda,\lambda^{-1})~|~\lambda\in \Lambda(\widetilde T)\}.
$$
Then we have
$$
  L_Z(\widetilde G)\cong \left(L(\widetilde G)\rtimes\Lambda(T)\right)/N\,.
$$
%

Choose a central extension
$\widehat\Lambda(T)$ of the lattice  $\Lambda(T)$ by $\C^*$. 
Any such central extension
is uniquely
determined by a skew-symmetric $\Z$-bilinear form (the commutator map) 
$\omega$ on
$\Lambda(T)$ which is defined by
$$
  \omega(\lambda,\mu) = 
    \hat\lambda\hat\mu\hat\lambda^{-1}\hat\mu^{-1}\,.
$$
Here, $\hat\lambda$ and  $\hat\mu$ are arbitrary lifts of
$\lambda,\mu\in\Lambda(T)$ to $\widehat\Lambda(T)$. Let   
$\widehat L^k(\widetilde G)$ be the central extension of $L(\widetilde G)$
of level $k$.
Suppose that $\widehat\Lambda(T)$ is a central extension of $\Lambda(T)$
such that its commutator map satisfies  
\begin{equation}\label{compat}
  \omega(\lambda,\mu)=(-1)^{k\langle\lambda,\mu\rangle}\quad~
   \text{ for all }~\lambda\in\Lambda(\widetilde T)\text{ and }
     \mu\in\Lambda(T)\,.
\end{equation}
Then one can construct a central
extension of $L_Z(\widetilde G)$ as follows: 
The group $\Lambda(T)\subset L_Z(\widetilde G)$ 
acts on $L(\widetilde G)$ by conjugation. This action uniquely 
lifts to an action of $\Lambda(T)$ on the central
extension $\widehat L^k(\widetilde G)$. We can consider the semidirect
product $\widehat L^k(\widetilde G)\rtimes\widehat\Lambda(T)$, 
where the action of
$\widehat\Lambda(T)$ on $\widehat L^k(G)$ 
factors through the action of $\Lambda(T)$. Now, the
lattice $\Lambda(\widetilde T)$ is a subgroup of $L(\widetilde G)$ so that the
restriction of the central extension of $L(\widetilde G)$ to 
this lattice yields
a central extension $\widehat\Lambda(\widetilde T)$ 
of $\Lambda(\widetilde T)$. On the other hand, we can restrict the central
extension $\widehat\Lambda(T)$ of $\Lambda(T)$ to the sublattice
$\Lambda(\widetilde T)$. 
The compatibility condition of equation (\ref{compat}) implies 
in particular that $\omega(\lambda,\mu)=(-1)^{k\langle\lambda,\mu\rangle}$ 
for all $\lambda,\mu\in\Lambda(\widetilde T)$. This implies 
that the two extensions of $\Lambda(\widetilde T)$ are equivalent 
(\cite{PS}, Proposition 4.8.1).
%
We may therefore consider the subgroup 
$$
  \widehat N=\{(\widehat\lambda,\widehat\lambda^{-1}) ~ | ~ 
 \widehat\lambda\in\widehat\Lambda(\widetilde T)\}\,
   \subset\,\widehat
    L^k(\widetilde G)\rtimes\widehat\Lambda(T)~.
$$
Now, using the full compatibility condition $(\ref{compat})$, 
one can check (\cite{T}, Proposition 3.3.1) 
that $\widehat N$ is a normal subgroup in 
$\widehat L^k(\widetilde G)\rtimes\widehat\Lambda(T)$ . 
Therefore, the quotient 
\begin{equation}\label{central}
  \widehat L^k_Z(\widetilde G) = 
    \left(\widehat L^k(\widetilde G)\rtimes\widehat\Lambda(T)\right)
     /\widehat N
\end{equation}
is a central extension of $L_Z(\widetilde G)$. 

We then have the
following theorem (\cite{T}, Theorem 3.2.1 and Proposition 3.3.1).
\begin{Theorem}\label{centralext}
  Every central extension of $L_Z(\widetilde G)$ is uniquely determined by the
  level $k$ of the corresponding central extension of $L(\widetilde G)$ and by
  a commutator map $\omega$ defining a central extension of 
  $\Lambda(T)$ which satisfies 
  the compatibility condition of equation
  (\ref{compat}). The corresponding central
  extension of $L_Z(\widetilde G)$ is the one described in equation
  (\ref{central}). 
\end{Theorem}
\begin{Remark}
  Note that Theorem \ref{centralext} restricts the possible levels at which
  central extensions of $L_Z(\widetilde G)$ can exist. For example,
  $L_{\langle-id\rangle}(SL(2,\C))$ does not 
  posses any central extensions of odd
  level. Indeed, we have $\Lambda(\widetilde T)=\alpha\Z$ and
  $\Lambda(T)=\frac\alpha2\Z$ with $\langle\alpha,\alpha\rangle=2$. Now, for
  odd level, the compatibility requirement of equation (\ref{compat}) 
  requires
  $\omega(\alpha,\frac\alpha2)=-1$ which is in contradiction to bilinearity
  and skew-symmetry of $\omega$.
\end{Remark}

\begin{Def}
  Let $k_f$ be the smallest level at which a central extension on 
  $L_Z(\widetilde G)$ exists. This $k_f$ is called the fundamental 
  level of $G$.
  Let $k_b$ be the smallest positive
  integer such that the restriction of $k_b\langle\cdot,\cdot\rangle$  to
  $\Lambda(T)$ is integral. $k_b$ is called the fundamental level of $G$.
\end{Def}

Obviously, for the fundamental level one has
$k_f\in\{1,2\}$. One can show (\cite{T}) 
that the basic level of $G$ is always a multiple of the fundamental 
level of $G$. The fundamental and basic levels of the simple Lie groups are
computed in \cite{T}. We list the basic levels in the appendix of this
paper. Finally, we have
(\cite{T}, Proposition 3.5.1):
\begin{Prop}
  A central extension of $L_Z(\widetilde G)$ is a pull-back of a central
  extension of $L(G)$ only if its level $k$ is a multiple of the basic level
  $k_b$ of $G$. Conversely, if $k_b|k$, the subgroup 
  $Z\subset\widehat L^k_Z(\widetilde G)$ 
  corresponding to the canonical embedding 
  $\widetilde G\hookrightarrow \widehat L^k_Z(\widetilde G)$ is central and 
  we have
  $$
    \widehat L^k_Z(\widetilde G)\cong
      \pi^*(\widehat L^k_Z(\widetilde G)/Z)~.
  $$
\end{Prop}    

\medskip

Let us fix a commutator map $\omega$ satisfying the compatibility
requirement from equation (\ref{compat}) for the rest of this paper

\medskip

From now on let us assume that $Z=\langle c\rangle$ is a cyclic group.
The group $\C$ acts naturally on $L_Z(\widetilde G)$ by
translations. This action factors through an action of
$\C/{\ord(c)}\Z$. We view  $\C/{\ord(c)}\Z$ as an $\ord(c)$-fold
covering of $\C/\Z\cong\C^*$, which we denote by 
$\widetilde{\C^*}$. 
Thus, we can  
define the semidirect product $L_Z(\widetilde G)\rtimes\widetilde{\C^*}$.
The action of $\widetilde{\C^*}$ on $L_Z(\widetilde G)$ described above 
lifts to any central extension of
$\widehat L^k_Z(\widetilde G)$ of $L_Z(\widetilde G)$ 
(\cite{T}, Proposition 3.4.1). 
So we can define the group 
$\widetilde L^k_Z(\widetilde G)
 = \widehat L^k_Z(\widetilde G)\rtimes\widetilde{\C^*}$.
Furthermore, if the level $k$ of the central extension 
$\widehat L^k_Z(\widetilde G)$ is a
multiple of the basic level $k_b$ of $G$, the group $Z=\langle c\rangle$ is a
central subgroup of $\widetilde L_Z(\widetilde G)$ and we can define 
$$
  \widetilde L^k(G)=\widetilde L^k_Z(\widetilde G)/Z~.
$$


Finally,  note that $\C^*$ acts naturally on the loop group
$L(G)$. However, contrary to the simply connected case, this action does not
necessarily lift to all central extensions of $L(G)$. In fact, we have
(\cite{T}, 3.5.10) 
\begin{Prop}
  The rotation action of $\C^*$ on $L(G)$ lifts to a central extension of
  $L(G)$ of level $k$ 
  if and only if $k\langle\lambda,\lambda\rangle\in2\Z$ for all 
  $\lambda\in\Lambda(T)$, i.e. if $\Lambda(T)$ endowed with
  $k\langle\cdot,\cdot\rangle$ is an even lattice.
\end{Prop}


\subsection{The adjoint action of $\widetilde L^k(G)$ on its Lie algebra}
\label{adjoint}
Suppose as before that $G=\widetilde G/Z$ where $\widetilde G$ 
is simply connected, and $Z= \langle c\rangle$ 
is a cyclic subgroup of the center 
of $\widetilde G$.
Consider the centrally extended loop group 
$\widetilde L^k_Z(\widetilde G) = 
 \widehat L^k_Z(\widetilde G)\rtimes\widetilde\C^*$ introduced in the
last section. Since the center of any Lie group acts trivially in the
adjoint representation, the adjoint action of 
$\widetilde L^k_Z(\widetilde G)$
on its Lie algebra $\widetilde L^k(\g)$ factors through
$L_Z(\widetilde G)\rtimes\widetilde{\C^*}$. The 
$\widetilde{\C^*}$-part acts by translations, so the only 
interesting part is the
action of $L_Z(\widetilde G)$.
Let $\zeta$ be an element of $L_Z(\widetilde G)$. Then the adjoint
action of $\zeta$ on $\widetilde L^k(\g)$ is given by 
(\cite{T} Corollary 3.4.2, \cite{PS})
\begin{multline}\label{ad}
  \Ad(\zeta): X+aC+bD\mapsto \zeta X\zeta^{-1}-\frac{b}{2\pi i}
  \dot\zeta\zeta^{-1} +  \\
  \big(a+\frac{k}{2\pi i}\int_0^{1}\
   \langle X(t),\zeta^{-1}(t)\dot\zeta(t)\rangle dt \\
  -\frac{kb}{8\pi^2 i^2}\int_0^{1}\langle\zeta^{-1}(t)\dot\zeta(t),
    \zeta^{-1}(t)\dot\zeta(t)\rangle dt \big)C + bD~.
\end{multline}
Here, $X$ is an element of the loop algebra $L(\g)$, and $\dot\zeta$
denotes the derivative of $\zeta$ with respect to $t$. Finally, as
before, we have identified $\C^*$ with $\C/\Z$ by identifying the
coordinate $z$ with $e^{2\pi i t}$.

We are interested in the action of a specific element of
$\sigma_c\in L_Z(\widetilde G)$ which is the product $\zeta_c w_c$ of an 
``open loop'' $\zeta_c$ and an element $w_c\in G$ which are defined as
follows. 
As before, let $\theta$ denote the
highest root of $\g$. The set of elements $\alpha_i\in\Pi$ which have
coefficient $m_j=1$ in the expansion $\theta=\sum_{j=1}^lm_j\alpha_j$
can be identified with the non-trivial elements of the center
of $\widetilde G$. Indeed, let $\{\lambda_j\}\subset\h$ be the dual
basis corresponding to $\Pi\subset\h^*$. Then the condition $m_j=1$
implies that $\exp(2\pi i\lambda_j)$ is an element of the center of
$\widetilde G$.   

Let $\alpha_c$ denote the root $\alpha$ which is
identified with the generator $c$ of $Z$ in this identification and let
$\lambda_c\in\h$ denote the corresponding fundamental weight of $\g$.  
There
exists a unique element $w_c\in W$ which permutes the set
$\Pi\cup\{-\theta\}$ and maps $-\theta$ to $\alpha_c$ 
(see \cite{T}, Proposition 4.1.2). Furthermore, 
let $\{e_\alpha~|~\alpha\in\Delta_+\}$ be a Chevalley basis
of the Borel subalgebra $\n\subset\g$. 
Then we can choose a representative $\bar w_c$ of 
$w_c$ in $N(T)$ such
that $\bar w_c(e_\alpha)=e_{w_c(\alpha)}$ for all $\alpha\in\Pi$. From
now on, we will denote both the Weyl group element $w_c$ as well as
its representative $\bar w_c\in N_G(T)$ simply by $w_c$. 
Finally, we can define an
``open loop'' $\zeta_c$ in  $G$ via 
$\zeta_c(t)=\exp(2\pi it\lambda_c)$.
Now, the element $\sigma_c$ is defined as
$$\sigma_c=\zeta_c^{-1}w_c.$$

The action of $\sigma_c$ on $\widetilde L^k(\g)$ can be described
explicitly in terms of the root space decomposition: To each root
$\alpha$ of $\g$ choose a co-root $h_\alpha\in\h$. Set
$h_{\alpha_0}=h_{-\theta}+kC$. Since $\langle\lambda_c,\alpha_c\rangle=1$
and $\langle\lambda_c,\alpha\rangle=0$ for all $\alpha\in\Pi$ with
$\alpha\neq\alpha_c$, we find that 
the ation of $\sigma_c$
on $\widetilde\h$ is
given by
$$D\mapsto D+\lambda_c-\frac{k}{2}\Vert\lambda_c\Vert^2C\,,$$
$$h_{\alpha_0}\mapsto h_{\alpha_c}\,,\quad
 h_{w_c^{-1}(-\theta)} \mapsto h_{\alpha_0}\,,\quad\text{and}$$
$$h_{\alpha_i}\mapsto h_{w_c(\alpha_i)}\quad
\text{for all other }i\,.$$
A set of generators for $\widetilde L(\g)_\pol$ is given by the set
$\{e_{\widetilde\alpha},f_{\widetilde\alpha}~|~\alpha\in\widetilde\Pi\}$ with
$e_{\widetilde\alpha_i}=e_{\alpha_i}\otimes 1$ for $1\leq i\leq l$ and
$e_{\widetilde\alpha_0}=e_{-\theta}\otimes z$, and accordingly
$f_{\widetilde\alpha_i}=f_{\alpha_i}\otimes 1$ for $1\leq i\leq l$ and
$f_{\widetilde\alpha_0}=f_{-\theta}\otimes z^{-1}$. 
It is straight forward
to check that $\sigma_c$ permutes the $e_{\widetilde\alpha_i}$ 
according to its
action on $\widetilde\Pi$ and similarly for the 
$f_{\widetilde\alpha_i}$.


\subsection{Integrable representations and characters}\label{reps}
As before, 
let $\theta$ denote the highest root of $\g$ and fix some non-negative
integer $k\in\Z_{\geq0}$. 
Let $P_+$ be the set of dominant weights of $\g$ with respect to
$\Pi$, and let $P_+^k$
denote the set of $\lambda\in P_+$ such that 
$\langle\lambda,\theta\rangle\leq k$. To each pair $(\lambda,k)$ with
$\lambda\in P^k_+$, we can associate an 
irreducible highest weight 
module $V_{\lambda,k}$
of $\widehat L(\g)_{\pol}$ 
such that the center of $\hat\g$ acts as the scalar $k$ (see \cite{K}). 
Letting $D$ act on the highest 
weight vector of  $V_{\lambda,k}$ 
as an arbitrary scalar uniquely determines an irreducible highest weight
representation of $\widetilde L(\g)_{\pol}$. 
We will denote by $V_{\lambda,k}$ the
highest weight representation of $\widetilde L(\g)_{\pol}$ such that 
$D$ acts trivially on the highest weight vector.

It was shown by H.~Garland \cite{G} that
$V_{\lambda,k}$ admits a positive definite
Hermitian form $(.,.)$ which is contravariant with respect to the anti-linear
Cartan involution on $\widetilde{L(\g)}_{\pol}$. Let us
denote by $V^{ss}_{\lambda,k}$ the $L^2$-completion of $V_{\lambda,k}$
with respect to this norm defined by the Hermitian form. That is, if 
$\{v_{\lambda k \mu i}\}_{i\in I(\lambda)}$ 
is an orthonormal basis of the weight subspace 
$V_{\lambda,k}[\mu]$ of $V_{\lambda,k}$, then 
$$V^{ss}_{\lambda,k} = \Big\lbrace\sum_{\mu,i} a_{\mu i}v_{\lambda k\mu i}~|~
\sum_{\mu,i} |a_{\mu i}|^2 < \infty \Big\rbrace.$$
Analogously, we define the analytic completion of $V_{\lambda,k}$ to be
the space
$$V^{an}_{\lambda,k} = \Big\lbrace\sum_{\mu,i} a_{\mu i}v_{\lambda k\mu i}~|~
\text{there exists a}~0<q<1
\text { s.t. } a_{\mu i}=O\left(q^{-D(\mu)}\right), \text{ }
\mu\to\infty \Big\rbrace,$$
where $D(\mu)$ denotes the (non-positive) degree of the weight $\mu$ in the 
homogeneous grading. 
By definition,  $V^{ss}_{\lambda,k}$ is a Hilbert space and 
$V^{an}_{\lambda,k}$ is a dense subspace in it. It is known 
(\cite{GW}, \cite{EFK}), that the action of $\widetilde L(\g)_{\pol}$ on 
$V_{\lambda,k}$ extends by continuity to an action of $\widetilde L(\g)$ on 
$V^{an}_{\lambda,k}$ , but not to an action on $V^{ss}_{\lambda,k}$.
 
\par

We now turn to the representation theory of the affine Lie groups. We
first consider the case that the corresponding finite dimensional Lie group is
simply connected. In this case, the following result is known (see
e.g. \cite{EFK}, Theorem 2.2 and Lemma 2.3).
\begin{Theorem}\label{trace}\quad
 \renewcommand{\theenumi}{\roman{enumi}}
\renewcommand{\labelenumi}{\textnormal{(\theenumi)}}
\begin{enumerate}
\item
  The action of the Lie algebra $\widehat L(\g)$ on $V^{an}_{\lambda,k}$
  uniquely integrates to an action of $\widehat L(G)$.
\item
  For any $q\in\C^*$ with $|q|<1$ and any $g\in\widehat L(G)$, the operator
  $gq^{-D}:V^{an}_{\lambda,k}\to V^{an}_{\lambda,k}$ uniquely 
  extends to a trace class operator on $V^{ss}_{\lambda,k}$.
\end{enumerate}
\end{Theorem}

Let $q\in\C^*$ and denote by $\widetilde L(G)_q$
the subset of $\widetilde L(G)$ of elements of the 
form $(g,q)$ with $g\in\widehat L(G)$.
This subset is invariant under conjugation in $\widetilde L(G)$.
Furthermore, let us introduce the semigroup 
$\widetilde L(G)_{<1}=\cup_{|q|<1}\widetilde L(G)_q$. As a manifold,
$\widetilde L(G)_{<1}$ is isomorphic to $\widehat L(G)\times D^*$,
where $D^*$ denotes the punctured unit disk in $\C$.
Since any element in $\widetilde L(G)_{<1}$ extends to a trace class operator
on $V^{ss}_{\lambda,k}$, we can introduce functions 
$\chi_{\lambda,k}:\widetilde L(G)_{<1}\to\C$ via
$$
  \chi_{\lambda,k}(g,q)=Tr|_{V^{an}_{\lambda,k}}(gq^{-D})\,.
$$
The function $\chi_{\lambda,k}$ is the character 
of the module $V_{\lambda,k}$.
\begin{Prop}[\cite{EFK}, Lemma 2.4, Proposition 2.5]\label{hol}
  The functions $\chi_{\lambda,k}$ are holomorphic and conjugacy invariant.
\end{Prop}

By definition, the central element $C$ of $\widetilde L(\g)$ acts 
on $V^{an}_{\lambda,k}$
by scalar multiplication with $k$. Therefore, if 
$\iota:\C^*\to\widehat L(G)$ denotes the identification of $\C^*$ with 
the center of $\widehat L(G)$, we have 
$\chi_{\lambda,k}(\iota(u)g) = u^k\chi_{\lambda,k}(g)$ for all $u\in \C^*$.

\medskip

Now let $G=\widetilde G/Z$ be non simply connected and let us fix once
and for all a central extension $\widehat L_Z(\widetilde G)$ of
$L_Z(\widetilde G)$ of level $k_f$.  
Suppose that $V$ is an irreducible representation of 
$\widehat L_Z(\widetilde G)$ of level $k$. Then $k$ is necessarily a
multiple of $k_f$. Restricting this representation to
the connected component of $\widehat L_Z(\widetilde G)$ containing
the identity  element yields a level $k$ representation of 
$\widehat L(\widetilde G)$. We shall always assume that we can
decompose $V$ as
\begin{equation}\label{decomp}
  V\cong\bigoplus_{\lambda\in I\subset P_+^k}V^{an}_{\lambda,k}
\end{equation}
as a representation of $\widehat L(\widetilde G)$. 
Such representations are called negative energy
representations of $\widehat L_Z(\widetilde G)$.

Remember the automorphism $\sigma_c$ introduced in the last
section. As we have seen, $\sigma_c$ acts as an automorphism of the
Dynkin diagram of $\widetilde L(\g)_\pol$. 
So, $\sigma_c$ permutes the set $P_+^k$
of level $k$-representations of $\widehat L(\g)_\pol$.  The following
Theorem is essentially Theorem 6.1 and Corollary
7.3. of \cite{T}
(note that \cite{T} deals with projective representations).


\begin{Theorem}
  Let $k$ be a multiple of the fundamental level of $G$. To each 
  $\sigma_c$-orbit $I\subset P_+^k$ there exist $\ord(c)/|I|$ different
  irreducible negative energy 
  representations of $\widehat L_Z(\widetilde G)$
  which, restricted to a representation of $\widehat L(\widetilde G)$
  decompose according to equation (\ref{decomp}).

  A representation $V_{I,k}$ of $\widehat L_Z(\widetilde G)$
  factors through a representation of 
  $\widehat L^{k_b}(G)$ if and only if 
  the level $k$ is a multiple of
  the basic level $k_b$ of $G$ and $Z$ acts trivially on $V_{I,k}$.
\end{Theorem} 

Suppose that the $\sigma_c$--orbit $I$ consists of a single
element $\lambda$. Then $\sigma_c$ acts on the highest weight subspace
$V_\lambda[\lambda]$. Since $V_\lambda[\lambda]$ is one--dimensional,
$\sigma_c$ acts by a scalar. From now on, we shall assume that
$\sigma_c$ acts as the identity on $V_\lambda[\lambda]$. This
determines the $\widehat L_Z(\widetilde G)$--module $V_{I,k}$
corresponding to $I$ and $k$ uniquely. 

Finally, letting $D$ act trivially on all highest weight vectors
$v_{\lambda,k}\in V_{I,k}$, we can extend the representation $V_{I,k}$
of $\widehat L_Z(\widetilde G)$ to a representation of 
$\widetilde L_Z(\widetilde G) = 
\widehat L_Z(\widetilde G)\rtimes\widetilde{\C^*}$ 
which factors through a representation of 
$\widetilde L(G)$ if $Z$ acts trivially on $V_{I,k}$ and $k_b|k$. 
(Recall that $\widetilde{\C^*}$ denotes the
$\ord(c)$-fold covering of $\C^*$ which acts on
$\widehat L_Z(\widetilde G)$ covering the natural $\widetilde{\C^*}$-action on
$L_Z(G)$). For any $\widetilde q\in\widetilde{\C^*}$, let us denote by
$q$ its image under the natural projection
$\widetilde{\C^*}\to \C^*$.

It is easy to see that $\sigma_c$ acts as a unitary operator in
$V_{I,k}$. So we have
\begin{Cor}
  Let $k$ be a multiple of $k_f$ (resp. $k_b$).
  For any $\widetilde q\in\widetilde{\C^*}$ with $|q|<1$ and any
  $g\in\widehat L_Z(\widetilde G)$ (resp. $g\in\widehat L(G)$),
  the operator
  $g\widetilde q^{-D}:V^{an}_{I,k}\to V^{an}_{I,k}$ uniquely 
  extends to a trace class operator on
  $V^{ss}_{I,k}=\bigoplus_{\lambda\in I}V^{ss}_{\lambda,k}$
\end{Cor}
As before, we can compute the traces of the trace class
operators. Abusing terminology slightly, we define the functions
$$
     \chi_{I,k}(g,\widetilde q^{-D}) = 
	Tr|_{V^{an}_{I,k}}(g\widetilde q^{-D})
$$
on  $\widetilde L_Z(\widetilde G)_{<1}$
resp. $\widetilde L(G)_{<1}$ not distinguishing the different spaces
they are defined on. (The semi groups
$\widetilde L_Z(\widetilde G)_{<1}$ and $\widetilde L(G)_{<1}$ are
defined analogously to the simply connected case). 
Finally, again using the fact that $\sigma_c$ acts as a unitary
operator on $V^{ss}_{I,k}$, we get
\begin{Cor}
 The functions $\chi_{I,k}$ are holomorphic and conjugacy invariant.
\end{Cor} 

The main goal of this paper is to compute the functions $\chi_{I,k}$
explicitly which will yield a generalization of the Kac-Weyl
character formula. This will be done in section \ref{Sec:twchar}. 
Here, we
state a trivial observation: The character $\chi_{I,k}$ restricted to the
component of $\widetilde L_Z(\widetilde G)_{<1}$ containing the
element $(e,q)$ is just the sum of the characters
$\chi_{I,k}=\sum_{\lambda\in I}\chi_{\lambda,k}$ of the characters
$\chi_{\lambda,k}$ of $\widetilde L(G)$ considered above. On the
other hand, $\sigma_c$ permutes the highest weight vectors of the
representations $V_{\lambda,k}$ with $\lambda\in I$. So if 
the $\sigma_c$-orbit $I$ consists of more
that one element, the function $\chi_{I,k}$
restricted to the  connected component of
$\widetilde L_Z(\widetilde G)_{<1}$ containing the element 
$(\sigma_c,\widetilde q)$ vanishes.



\section{Conjugacy classes}\label{conjugacy}
\subsection{Conjugacy classes and principal bundles}
Since the characters $\chi_{I,k}$ are conjugacy invariant functions on
$\widetilde L^k(G)_q$, it is necessary to have a good understanding of the
conjugacy classes in $\widetilde L^k(G)_q$ in order to understand the
characters. The fundamental result in this direction is an observation due to 
E.~Looijenga which gives a one-to-one correspondence between the 
$L(G)$-conjugacy
classes in $L(G)\times\{q\}\subset L(G)\rtimes\C^*$ and the isomorphism
classes of holomorphic 
principal $G$-bundles over the elliptic curve $E_q=\C^*/q^\Z$. To be more
precise, let $G=\widetilde G/Z$, as above. Then, up to
$C^\infty$-isomorphism,
every principal $G$-bundle over $E_q$ is determined by its topological class,
which is an element in $\pi_1(G)\cong Z$. We can classify holomorphic
principal $G$-bundles of a fixed topological class $c$ as follows. 
Consider the connected component $L(G)_{c}$ of the loop group 
$L(G)$ which corresponds to $c$. 
The group $L(G)\rtimes\C^*$ acts on the set
$L(G)_{c}\times\{q\}\subset L(G)\rtimes\C^*$ by conjugation.
Looijenga's observation gives a one-to-one
correspondence between the set of holomorphic $G$-bundles on $E_q$ of
topological type $c$ and the set of $L(G)\times\{1\}$-orbits in
$L(G)_{c}\times\{q\}$. This correspondence comes about as follows. For any
element 
$(g,q)\in L(G)_{c}\times\{q\}$ consider the $G$-bundle $B_g$ over $E_q$
which is defined as follows. View $E_q$ as the annulus $\vert q\vert\leq\vert
z\vert\leq 1$ in the complex plane with the boundaries identified via
$z\mapsto qz$. Then take the trivial bundle over the annulus and define the
bundle $B_g$ over $E_q$ by describing the gluing map which identifies the
fibers over the points identified under $z\mapsto qz$. This is given by
$f(qz)=g(z)f(z)$, where $f$ takes values in $G$. Obviously, this construction
gives  a holomorphic $G$-bundle of topological type $c$ and the following
theorem which is due to E. Loojienga 
is not hard to prove (see e.g. \cite{EF}).
\begin{Theorem}\label{bundles}\quad
\renewcommand{\theenumi}{\roman{enumi}}
\renewcommand{\labelenumi}{\textnormal{(\theenumi)}}
\begin{enumerate}
\item
  Two elements $(g_1,q),\,(g_2,q)\in L(G)_c\times\{q\}$ are conjugate under
  $L(G)\times\{1\}$ if and only if the corresponding holomorphic
  $G$-bundles $B_{g_1}$ and
  $B_{g_2}$ are isomorphic.  
\item
  For any holomorphic 
  $G$-bundle over $E_q$ of topological type $c$, there
  exists an element $g\in L(G)_c$ such that $B\cong B_g$.
\end{enumerate}
\end{Theorem}

Following \cite{EFK}, we call an $L(G)\times\{1\}$-orbit 
in $L(G)_c\times\{q\}$
semisimple 
if the corresponding principal $G$-bundle over $E_q$ 
comes from a representation of the fundamental group of $E_q$
inside a maximal compact subgroup of $G$. We call an element
$(g,q)\in L(G)\times\{q\}$ semisimple if the corresponding 
$L(G)\times\{1\}$-orbit is semisimple.
It is known that almost every conjugacy class is semisimple. To be more
precise, 
if $\{\mathcal{B}_t\}_{t\in\mathcal{T}}$ is a holomorphic family of holomorphic
principal bundles on $E_q$ 
parametrized by a complex space $\mathcal{T}$, then the set
subset $\mathcal{T}_0$ of bundles which are flat and unitary is nonempty and
Zariski open in $\mathcal{T}$ (see \cite{R}).

\subsection{The simply connected case}
The set of semisimple $L(G)\times\{1\}$-orbits in $L(G)_c\times\{q\}$ 
can be described more explicitly. We
start with the case that $G$ is simply connected (the non-simply
connected case will be considered in section \ref{nsc}).
The affine Weyl group
$\widetilde W=W\ltimes\Lambda(T)$ can be identified with the group
$N_{L(G)\rtimes\C^*}(T\times\C^*)/(T\times\C^*)$ where
$N_{L(G)\rtimes\C^*}(T\times\C^*)$ denotes the normalizer of
$T\times\C^*$ in $L(G)\rtimes\C^*$ (see \cite{PS}). 
In this sense, 
$\widetilde W$ is the Weyl group of $\widetilde L(G)$. 
It acts on the torus $T\times\C^*$ via
$(w,\beta):(\xi,q)\mapsto (w(\xi)q^{-\beta},q)$. The
following proposition follows from the definition of semisimplicity and the
classification of stable and unitary $G$-bundles over $E_q$ (\cite{EF},
\cite{FM1}). 
\begin{Prop}\label{conclasses}
  Let $G$ be simply connected. 
  Every semisimple element in $L(G)\times\{q\}$ 
  is $L(G)\times\{1\}$--conjugate to an element of the form $(\xi,q)$
  with $\xi\in T$. Two elements $(\xi_1,q)$ and  $(\xi_2,q)$ with
  $\xi_1,\xi_2\in T$ are conjugate
  if and only if they are in the same orbit under the action of
  the affine Weyl group $\widetilde W$ on $T\times\C^*$.
\end{Prop}
\begin{Remark}\label{mod}
  We have $T=(\C^{*})^r$ for some $r\in\N$. The set $T/q^{\Lambda(T)}$
  is an Abelian variety isomorphic to the product
  $E_q\otimes_\Z\Lambda(T)$ and 
  there is a natural action of $W$ on  $E_q\otimes_\Z\Lambda(T)$. 
  In this way, the
  set of semisimple $L(G)\times\{1\}$-conjugacy classes in $L(G)\times\{q\}$ 
  can be identified with the set $E_q\otimes_\Z\Lambda(T)/W$.
\end{Remark}

Conjugacy classes in the centrally extended group $\widetilde L^k(G)$
can be described as follows.
We are interested in the semisimple conjugacy classes in
$\widetilde L^k(G)_q$, i.e. conjugacy classes which project to semisimple
conjugacy classes under the natural projection  
$\widetilde L^k(G)_q\to L(G)\times\{q\}$. The set $\widetilde L^k(G)_q$ is the
total space of a fiber-bundle over $L(G)\times\{q\}$ with fiber $\C^*$ and
conjugation with an element of $L(G)\times\{1\}$ induces an automorphism of
this bundle. So the set of semisimple $L(G)\times\{1\}$-orbits in 
$\widetilde L^k(G)_q$ 
will be the total space of a $\C^*$--bundle over the set of semisimple
$L(G)\times\{1\}$-orbits in $L(G)_q$. Using equation (\ref{ad}), one can
describe the bundle explicitly. 
\par
Let $\widehat T_q$ be the set of all elements in $\widetilde L^k(G)_q$ which
project  
to an element of the form $(\xi,q)$ with $\xi\in T$. As a complex manifold
this set is isomorphic to $T\times\C^*\times\{q\}$. 
\begin{Prop}\label{conprop}
  The set of semisimple 
  conjugacy classes in $\widetilde L^k(G)_q$ is given by the quotient 
  $\widehat T_q/\widetilde W$, where 
  $(w,\beta)\in \widetilde W=W\ltimes\Lambda(T)$ acts on 
  $\widehat T_q$ as follows:
  $$
    (w,1)(\xi,u,q) = (w(\xi),u,q),
  $$
  $$
    (1,\beta)(\xi,u,q) = \left(\xi q^{-\beta},
      uq^{-(k/2)\langle\beta,\beta\rangle}\beta(\xi^{k}),q\right)\,.
  $$
  Here, $\beta(\cdot)$ denotes the value of $\beta$ as a character of $T$ and 
  $\h$ and $\h^*$ are identified via the normalized invariant bilinear
  form on $\g$.  
\end{Prop}
\begin{Remark}\label{Lk}
  The quotient $\widehat T_q/\Lambda(T)$ described in Proposition \ref{conprop}
  is the set of non-zero vectors in 
  a holomorphic 
  line bundle $\L^k$ over the Abelian variety $E_q\otimes_\Z\Lambda(T)$. 
  The action of the finite Weyl group $W$ on this variety induces 
  an action of $W$ on $\L^k$.
\end{Remark}

\subsection{The non-simply connected
  case}\label{nsc}
If $G$ is not simply connected, the set of semisimple 
conjugacy classes in the connected
component $L(G)_c\times\{q\}$ of $L(G)\rtimes\C^*$ can be described
similarly to the simply connected case although the corresponding
analysis is slightly more involved. Throughout this section 
let us assume that $G$
is of the form $G=\widetilde G/Z$, where $\widetilde G$
is simply connected and $Z=\langle c\rangle$ 
is a cyclic subgroup of the center of 
$\widetilde G$.

Recall the the definition of the element $\sigma_c=\zeta_c^{-1}w_c$
from section \ref{adjoint}. The element $(\sigma_c,1)\in
L(G)\rtimes\C^*$ acts on the torus $T\times\C^*\subset
L(G)\rtimes\C^*$ by conjugation. We denote by $\big(T
\times\C^*\big)^{\sigma_c}_0$ the connected component of the fixed
point set $\big(T
\times\C^*\big)^{\sigma_c}$ which contains the identity element. 
%

Denote by
$L(G)_c\times\{q\}$ the connected component of
$L(G)\times\{q\}$ containing the element $(\sigma_c,q)$. 
The following lemma follows from the definition of semisimplicity and
the classification of flat and unitary $G$-bundles of topological type
$c$ on the elliptic curve $E_q$ (\cite{S} Equation 2.5, \cite{FM1}).

\begin{Lemma}
  Every semisimple element in $L(G)_c\times\{q\}$ is conjugate
  under $L(G)$ to an element of the form 
  $(\sigma_c\xi,q)\in(\sigma_c,1)\big(T\times\C^*\big)^{\sigma_c}_0$.
\end{Lemma}

It remains to check, which elements of
$(\sigma_c,1)\big(T\times\C^*\big)^{\sigma_c}_0$ are conjugate in
$L(G)\times\C^*$. Obviously, it is enough, to consider conjugation
with elements of $L(G)_0\rtimes\C^*$, where $L(G)_0$ denotes the
connected component of $L(G)$ containing the identity. In particular,
$L(G)_0$ consists of loops which are contractible.
So we have to study the ``twisted Weyl group''
$$
  \widetilde W_{\sigma_c} =
  N_{L(G)_0\rtimes\C^*}\left((\sigma_c,1)\big(T\times\C^*\big)^{\sigma_c}_0
  \right) / \big(T\times\C^*\big)^{\sigma_c}_0 
$$
and its action on the set $(\sigma_c,1)\big(T\times\C^*\big)^{\sigma_c}_0$.
First, note that if some $(g,q)$ normalizes
$(\sigma_c,1)\big(T\times\C^*\big)^{\sigma_c}_0$, then it normalizes
$\big(T\times\C^*\big)^{\sigma_c}_0$.

\begin{Lemma}\label{reg}
  If some element $(g,q)\in L(G)\times\C^*$ normalizes
  the torus $\big(T\times\C^*\big)^{\sigma_c}_0$, then it has to normalize
  the torus $T\times\C^*$ as well.
\end{Lemma}
\begin{proof}
  Conjugation by an element $(g,q)\in L(G)\rtimes\C^*$ 
  induces an automorphism of the Lie algebra
  $L(\g)\oplus \C D$ of $L(G)\rtimes\C^*$ which we will denote by $g_q$. 
  Consider the root space decomposition 
  $$
     L(\g)_\pol\oplus \C D=(\h\oplus\C D)
       \oplus\bigoplus_{\alpha\in\widetilde\Delta} L(\g)_\alpha\,.
  $$
  The automorphism $\sigma_c$ acts on $L(\g)\oplus \C D$ leaving
  $\h\oplus\C D$ invariant.
  We can choose an element $X\in(\h\oplus\C D)^{\sigma_c}$ 
  such that $\alpha(X)\neq 0$ for
  all $\alpha\in\widetilde\Delta$. This choice of $X$ insures that the
  condition $[X,Y]=0$ already
  implies $Y\in\h\oplus\C D$. But for $Y\in\h\oplus\C D$ and 
  $g\in N_{L(G)\rtimes\C^*}(T\times\C^*)^{\sigma_c}_0$
   we have
  $[X,g_q(Y)]=0$ so that $g_q$ indeed leaves $\h\oplus\C D$ invariant.
  Finally, we can use the exponential map 
  $\exp:L(\g)\oplus \C D\to L(G)\rtimes\C^*$ to go back to the group level.
%
%
\end{proof}
Let us denote by $\widetilde W$ the group 
$$
  \widetilde W =   N_{L(G)_0\rtimes\C^*}(T\times\C^*)
  / (T\times\C^*)\,.
$$
Then we have the following Lemma. 
\begin{Lemma}
  The twisted Weyl 
  group $\widetilde W_{\sigma_c}$ is isomorphic to the semidirect
  product
  $$ 
    \widetilde W_{\sigma_c}\cong \widetilde
    W^{\sigma_c}\ltimes\big((T\times\C^*) /
    (T\times\C^*)^{\sigma_c}_0\big)^{\sigma_c}\,. 
  $$
  of the $\sigma_c$-invariant
  part  $\widetilde W^{\sigma_c}$ of $\widetilde W$ and  
  the finite group
  $\big( (T\times\C^*) / (T\times\C^*)^{\sigma_c}_0\big)^{\sigma_c}$.
\end{Lemma}
\begin{proof}
  Lemma \ref{reg} allows to define a map
  $$ 
    \varphi: \widetilde W_{\sigma_c} \to \widetilde W
  $$
  via
  $$
    (g,q)(T\times\C^*)^{\sigma_c}_0\mapsto (g,q)(T\times\C^*)\,.
  $$
  It is easy to check that 
  $$
   (\sigma_c,1)\varphi\left((g,q)(T\times\C^*)^{\sigma_c}_0\right)
    (\sigma_c,1)^{-1} = 
     \varphi\left((g,q)(T\times\C^*)^{\sigma_c}_0\right) 
  $$
  so that $\varphi$ defines a map 
  $\widetilde W_{\sigma_c} \to \widetilde W^{\sigma_c}$. 
  The kernel of $\varphi$ is given by 
  $$
    ker(\phi)=\big((T\times\C^*)/(T\times\C^*)^{\sigma_c}_0\big)^{\sigma_c}\,.
  $$ 
  Finally, one can show that
  that the map
  $\varphi:\widetilde W_{\sigma_c} \to \widetilde W^{\sigma_c}$ is
  surjective and that the exact sequence
  $$
  \{1\}\longrightarrow \big((T\times\C^*) /
    (T\times\C^*)^{\sigma_c}_0\big)^{\sigma_c} \longrightarrow 
    \widetilde W_{\sigma_c}\overset{\varphi}{\longrightarrow}
       \widetilde W^{\sigma_c}\longrightarrow\{1\}
  $$
  splits (see e.g. \cite{W} or \cite{M} for the case of finite Weyl groups).
\end{proof}

It remains to describe the action of the groups $\widetilde W^{\sigma_c}$
and $\left( (T\times\C^*) / (T\times\C^*)^{\sigma_c}_0\right)^{\sigma_c}$
on $(\sigma_c,1)(T\times\C^*)^{\sigma_c}_0$. 
Note that conjugation by $(\sigma_c,1)$
maps an element $(\xi,q)\in T\times\C^*$ to $(w_c(\xi)q^{-\lambda_c},q)$. 
Let us choose some $h_0\in\h_\R$ such that $w_c(h_0)=h_0+\lambda_c$. Then 
we can define a bijective map 
$$T^{w_c}_0\times\C^*\to\big(T\times\C^*)^{\sigma_c}_0$$ 
via 
$$(\xi,q)\mapsto(\xi q^{h_0},q)\,.$$
This shows that 
$\big( (T\times\C^*) / (T\times\C^*)^{\sigma_c}_0\big)^{\sigma_c}$ is
in fact isomorphic to $\big(T/T^{w_c}_0\big)^{w_c}$. 
Let $T_{w_c}=T/(id-w_c)T$ denote the torus of co-invariants under $w_c$. There
is a natural projection $T^{w_c}_0\to T_{w_c}$ whose kernel is isomorphic
to $\big(T/T^{w_c}_0\big)^{w_c}$. This projection allows us to embed
the co-character lattice $\Lambda(T_{w_c})$ of $T_{w_c}$
into $\h^{w_c}$ such that 
$\Lambda(T_0^{w_c})\subset\Lambda(T_{w_c})$.
The group $\big(T/T^{w_c}_0\big)^{w_c}$
acts on $(\sigma_c,1)\big(T\times\C^*)^{\sigma_c}_0$
via translations.
The action of $\widetilde W^{\sigma_c}$ on
$(\sigma_c,1)\big(T\times\C^*)^{\sigma_c}_0$ is described in the
following lemma.
\begin{Lemma}\label{tildeW}
  Suppose that $\sigma_c$ is not the order $n$ automorphism
  of the extended Dynkin diagram corresponding to $A_{n-1}$.
  Then the group $\widetilde W^{\sigma_c}$ is isomorphic to a semidirect
  product
  $$
    \widetilde W^{\sigma_c}\cong W_0\ltimes \Lambda(T_{w_c})\,.
  $$
  Here, $W_0$ is a finite Weyl group acting irreducibly on $\h^{w_c}$.
  The action of $\widetilde W^{\sigma_c}$  
   on $(\sigma_c,1)(T\times\C^*)^{\sigma_c}_0$ is given by
  $$
    (w,1)\left((\sigma_c,1)(\xi q^{h_0},q)\right) =
      (\sigma_c,1)(w(\xi) q^{h_0},q)
  $$
  and  
  $$
    (1,\beta)\left((\sigma_c,1)(\xi q^{h_0},q)\right) = 
      (\sigma_c,1)(\xi q^{-\beta}q^{h_0},q)\,.
  $$
\end{Lemma}
\begin{proof}
  Fix some $q\neq 1$ in $\C^*$ and view the tangent space of
  $T\times\{q\}$ as an affine subspace of the tangent space $\h\oplus\C D$ of
  $T\times\C^*$. We can identify this space with $\h$.
  Since the Weyl group $\widetilde W$ maps  $T\times\{q\}$ to itself, it
  acts on $\h$ by affine transformations. It is a standard fact that
  $\widetilde W$ acts on $\h$ by affine
  reflections, and its lattice of translations is given by
  $\Lambda(\widetilde T)$. Let $\mathfrak{a}$ be a fundamental domain
  for the induced action of $\widetilde W$ on $\h_\R$. Then
  $\widetilde W$ is generated by the reflections in the walls of
  $\mathfrak{a}$. Now, $\sigma_c$ induces an affine map on $\h_\R$ which
  maps the set of reflection hyperplanes of $\widetilde W$ to
  itself. Furthermore, we can choose a fundamental domain
  $\mathfrak{a}_0$ which is mapped to itself by $\sigma_c$. So
  $\mathfrak{a}_0\cap\h^{\sigma_c}_\R\neq\emptyset$. We claim that the
  action of $\widetilde W^{\sigma_c}$ on $\h_\R^{\sigma_c}$ is generated
  by the set of reflections in the walls of
  $\mathfrak{a}_0\cap\h^{\sigma_c}_\R$. Indeed, let $\h_\nu$ be a wall of
  $\mathfrak{a}_0\cap\h^{\sigma_c}_\R$ 
  which is the intersection of the walls $\h_1,\ldots,\h_r$
  of $\mathfrak{a}_0$. 
  Since each $x\in\h_1\cap\ldots\cap\h_r$ is fixed under $\sigma_c$, the map
  $\sigma_c$ permutes the hyperplanes $\h_1,\ldots,\h_r$. Furthermore, since
  $\h_1\cap\ldots\cap\h_r$ is supposed to be a wall of 
  $\mathfrak{a}_0\cap\h^{\sigma_c}_\R$
  (i.e. an affine  subspace of $\h^{\sigma_c}$ of codimension $1$), we can
  assume the $\h_i$ to consist of a single $\sigma_c$-orbit. Let
  $s_i:\h_\R\to\h_\R$ denote 
  the reflection in the affine hyperplane $\h_i$. Since $\sigma_c$ is
  not the order $n$-automorphism of the extended Dynkin diagram
  corresponding to $A_{n-1}$, we have that either for all simple roots
  $\widetilde\alpha_i\in\widetilde\Pi$, the simple roots $\widetilde\alpha_i$
  and $\sigma_c(\widetilde\alpha_i)$ are non connected in the Dynkin
  diagram of $\widetilde\Delta$ in which case $s_i$ and
  $s_{\sigma_c(i)}$ commute.  Or the $\sigma_c$-orbit through 
  $\widetilde\alpha_i$ consists of exactly two elements $s_1$ and
  $s_2$ say, and $(s_1s_2)^3=1$. 
  In the first case, $s_1s_2\cdots s_r$ commutes
  with $\sigma_c$. Furthermore, the restriction of $s_1s_2\cdots s_r$ leaves
  an affine hyperplane of $\h^{\sigma_c}$ invariant so that it acts as an
  affine reflection on $\h^{\sigma_c}$. In the second case, we have to
  consider the element $s_1s_2s_1\in\widetilde W$ which commutes with
  $\sigma_c$ and acts as an affine reflection on $\h^{\sigma_c}$.  
  Since
  $\mathfrak{a}_0\cap\h^{\sigma_c}_\R$  
  is a fundamental domain for the action of
  $\widetilde W^{\sigma_c}$ acting on $\h^{\sigma_c}$, the subgroup of
  $\widetilde W^{\sigma_c}$  generated
  by the reflections  in the walls of $\mathfrak{a}_0\cap\h^{\sigma_c}_\R$
  generate the whole of  $\widetilde W^{\sigma_c}$. 

  Finally, it is straight
  forward to check that the lattice of translations of the action of
  $\widetilde W^{\sigma_c}$ on $\h^{\sigma_c}$ is given by the lattice 
  $$
    \Lambda(T_{w_c}) = \lbrace\frac{1}{\ord(w_c)}\sum_{i=1}^{\ord(w_c)}
      w_c^i(\beta)~|~\beta\in\Lambda(\widetilde T)\rbrace\,.
  $$
\end{proof}
Putting everything together, we get the analogue of Proposition 
\ref{conclasses}
in the non-simply connected case.
\begin{Prop}\label{conclasses2}
  Let $G$ be of the form $G=\widetilde G/Z$.
  Every semisimple element in $L(G)_c\times\{q\}$ is 
  conjugate to an element of the form $(\sigma_c\xi q^{h_0},q)$,
  where $\xi\in T^{w_c}_0$ and 
  $h_0\in\h_\R$ is chosen such that $w_c(h_0)=h_0+\lambda_c$.

  Two elements $(\sigma_c\xi_1q^{h_0},q)$ and  
  $(\sigma_c\xi_2q^{h_0},q)$ are conjugate
  if and only if they are in the same orbit under the action of the group
  $\widetilde W_{\sigma_c}$ 
  on $(\sigma_c,1)(T\times\C^*)^{\sigma_c}$.
\end{Prop}

\begin{Remark}\label{mod2}
  Remember that we have identified $\big(T/T^{w_c}_0\big)^{w_c}$ with a finite
  subgroup of $T^{w_c}_0$. The quotient
  $T^{w_c}_0/\big(T/T^{w_c}_0\big)^{w_c}$ is isomorphic to the torus $T_{w_c}$
  of co-invariants. So analogously to Remark \ref{mod},  the
  set of semisimple conjugacy classes in $L(G)_c\times\{q\}$ 
  can be identified with the quotient $E_q\otimes_\Z\Lambda(T_{w_c})/W_0$.
\end{Remark}

Finally, let us consider semisimple conjugacy classes in a central
extension $\widetilde L^k(G)$. Remember that 
$\widetilde L^k(G)$ is given as a semidirect product $\widetilde
L^k(G)\cong\widehat L^k(G)\rtimes\widetilde{\C^*}$ where
$\widehat L^k(G)$  denotes a level $k$ central extension of $L(G)$
corresponding to some commutator map $\omega$,  and
$\widetilde{\C^*}$ is the $\ord(c)$-fold covering of $\C^*$.
Let $L(G)_c$ denote the connected component of the loop group $L(G)$
corresponding to the element $c\in Z$ and denote by $\widehat L^k(G)_c$ the
set of elements in $\widehat L^k(G)$ which project to an element in
$L(G)_c$.
Furthermore, fix some $\widetilde q\in\widetilde{\C^*}$. Since
conjugation in $\widetilde L^k(G)$ leaves the sets $\widehat
L^k(G)\times\{\widetilde q\}$ invariant, we can consider
the connected component $\widehat L^k(G)_c\times\{\widetilde q\}$
which is also invariant under conjugation.
The set of semisimple conjugacy classes in this component 
will be the total
space of a $\C^*$-bundle over the space of semisimple conjugacy classes in
$L(G)_c\times\{\widetilde q\}$. For any element $\widetilde
q\in\widetilde{\C^*}$, let $q\in\C^*$ denote the image of
$\widetilde q$ under the natural projection $\widetilde{\C^*}\to\C^*$.
Since $\widetilde q$ acts on $L(G)$ via rotation by $q$, the set of
semisimple conjugacy classes in $L(G)_c\times\{\widetilde q\}$ is
isomorphic to the set of semisimple conjugacy classes in
$L(G)_c\times\{q\}$.

Recall from Lemma \ref{tildeW} and Remark \ref{mod2} 
that 
$\widetilde W_{\sigma_c} = 
  \left(W_0\ltimes \Lambda(T_{w_c})\right)\ltimes(T/T^{w_c}_0)$, where
$\Lambda(T_{w_c})$ is identified with a sublattice of $\h_\R^{w_c}$.
 
Let us denote by $\widehat T^{w_c}_q$ the set of elements of $\widehat
L^k(G)_c\times\{\widetilde q\}$ 
which project to an element of the form 
$(\sigma_c\xi q^{h_0},\widetilde q)$
with $\xi\in T^{w_c}_0$. As a complex manifold, this set is isomorphic to
$T^{w_c}_0\times\C^*$. 
Using Lemma \ref{tildeW}, Remark \ref{mod2} and formula (\ref{ad}), we
find the analogue of Proposition \ref{conprop} in the non-simply
connected case.
\begin{Prop}\label{conprop2}
  The set of semisimple 
  $L(G)$-conjugacy classes in the connected component of 
  $\widehat L^k(G)\times\{\widetilde q\}$ 
  corresponding to $c\in Z$ is given by the quotient 
  $\widehat T^{w_c}_q/\widetilde W_{\sigma_c}$, where 
  $(w,\xi_0,\beta)\in \widetilde W_{\sigma_c} = 
    \left(W_0\ltimes \Lambda(T_{w_c})\right)\ltimes(T/T^{w_c}_0)$
   acts as
  follows:
  $$
    (w,1,1)(\xi,u,\widetilde q) = (w\xi,u,\widetilde q)\,,
  $$
  $$
    (1,\beta,1)(\xi,u,\widetilde q) = \left(\xi q^{-\beta},
    uq^{-(k/2)\langle\beta,\beta\rangle}\beta(\xi^{k}),\widetilde q\right)
    \quad\text{for } \beta\in\Lambda(T^{w_c}_0)\,.
  $$
  and
  $$
    (1,1,\xi_0)(\xi,u,\widetilde q)=(\xi\xi_0,u,\widetilde q)\,,
  $$
  Here, $\beta(\cdot)$ denotes the value of $\beta$ 
  as a character of $T_{w_c}$ 
  and 
  $\h$ and $\h^*$ are identified via the normalized 
  invariant bilinear form on $\g$.
\end{Prop}
\begin{Remark}\label{Lk2}
  The quotient $\widehat T^{w_c}_q/\Lambda(T_{w_c})$ 
  described in Proposition \ref{conprop2}
  is the set of non-zero vectors in 
  a line bundle $\L_{w_c}^k$ over the Abelian 
  variety $E_q\otimes_\Z\Lambda(T_{w_c})$
  introduced in Remark \ref{mod}. 
  The action of the finite Weyl group $W_0$ on this variety induces 
  an action of $W_0$ on $\L_{w_c}^k$.
\end{Remark}



\section{The differential equation for affine characters}
\label{proof}

\subsection{The characters as sections of a line bundle}\label{charsec}
Throughout this
section, we shall assume that $G$ is of the form 
$G=\widetilde G/Z$, where $Z=\langle c\rangle$ is a cyclic subgroup of
the center of $\widetilde G$. The simply connected case follows
from the calculations below by setting $c=id$.
As before, let $k$ be a multiple of the basic level of $G$, and let
$\chi_{I,k}$ denote the character of the 
$\widetilde L^{k_b}(G)$-module $V_{I,k}$. Here, $I$ is a $\sigma_c$-orbit in
$P^k_+$. So the representation $V_{I,k}$ restricted to the connected
component $\widetilde L^{k_b}(G)_0$ of $\widetilde L^{k_b}(G)$
containing the identity decomposes into a direct sum
$V_{I,k}=\oplus_{\lambda\in I}V_{\lambda,k}$ of irreducible highest
weight representations of $\widetilde L^{k_b}(G)_0$. 
We have seen in the end of section \ref{reps} that if $I$
consist of more than one element, the character $\chi_{I,k}$
restricted to the connected component of $\widetilde L^{k_b}(G)$
corresponding to the element $c$   
vanishes. So the only interesting characters on this component are the
ones coming from 
representations $V_{I,k}$, where $I$ consists of a single element
$\lambda$ which is necessarily invariant under $\sigma_c$. Let us
denote the corresponding character by $\chi_{\lambda,k}$.

As before, let $L^k(G)_c$ denote the connected component of $L(G)$
containing the element $\sigma_c$, and let $\widehat L^{k_b}(G)_c$ be the
connected component of $\widehat L^{k_b}(G)$ which consists of elements
which project to some $g\in L(G)_c$. Fix some 
$\widetilde q\in\widetilde\C^*$ and
let $q$ denote the image of $\widetilde q$ under the natural
projection $\widetilde\C^*\to\C^*$. Let us assume that $|q|<1$. Then
the character $\chi_{\lambda,k}$ defines a holomorphic function on
$\widehat L^{k_b}(G)\times\{\widetilde q\}$. Since almost every element
in $\widehat L^{k_b}(G)\times\{\widetilde q\}$ is semisimple,
the function $\chi_{\lambda,k}$ is uniquely determined by its values on 
$\widehat T^{w_c}_q\times\{\widetilde q\}$. 
As before, fix some $h_0\in\h$ such that $w_c(h_0)=h_0+\lambda_c$.
Remember that $\widehat T^{w_c}_q\times\{\widetilde q\}$ 
was defined as the set of elements
of $\widehat L^{k_b}(G)_c\times\{\widetilde q\}$ 
which project to an element of the form
$(\sigma_c\xi q^{h_0},\widetilde q)$ with $\xi\in T^{w_c}_0$
under the natural projection
$\widehat L^{k_b}(G)\times\{\widetilde q\} \to 
 L(G)\times\{\widetilde q\}$. 

Finally, recall the identification $\iota:\C^*\to \widehat L^{k_b}(G)$
of $\C^*$ with the center of 
$\widehat L^{k_b}(G)$. 
If $g=(\sigma_c\xi q^{h_0},u,\widetilde q)\in\widehat T^{w_c}_q$, 
we can use the identity 
$$
  \chi_{\lambda,k}(\sigma_c\xi q^{h_0},u,\widetilde q) = 
  u^k\chi_{\lambda,k}(\sigma_c\xi q^{h_0},1,\widetilde q)
$$
to get rid of the central variable. Thus, for fixed $\widetilde q$, 
we can view the character
$\chi_{\lambda,k}$ as a section of the  line bundle
$\L^k_{w_c}$ introduced in Remark \ref{Lk2}.  
We shall abuse notation slightly and denote
this section by $\chi_{\lambda,k}$ as well. So locally, we can view
$\chi_{\lambda,k}$ as a function on
$\big(T\times\C^*\big)^{\sigma_c}_0$.
We will derive a differential equation for this function in
section \ref{theDE} by varying the variable $\widetilde q$.

Let us change our notation slightly. For $q\in\C^*$ with $|q|<1$ 
fix $\tau\in\C$ with $Im(\tau)>0$ such that
$q=e^{2\pi i\tau}$. Let $L\subset\C$ be the lattice generated by $1$ and 
$\tau$. Then the elliptic curve $E_\tau=\C/L$
is isomorphic to the curve $E_q$ considered in the last paragraph via the 
map $x\mapsto e^{2\pi i x}$. This allows to identify the Abelian
variety  $E_q\otimes_\Z\Lambda(T_{w_c})$ with
$\h^{w_c}/\left(\Lambda(T_{w_c})\oplus\tau\Lambda(T_{w_c})\right)$. In this
identification, we can view the character $\chi_{\lambda,k}$ at some
fixed $\widetilde q$ as a function  
$\chi^{\sigma_c}_{\lambda,k}(\cdot;\widetilde q)$
on $\h^{w_c}$ which has to satisfy
the following identity:
$$
  \chi^{\sigma_c}_{\lambda,k}(h+\beta+\tau\beta';\widetilde q)=
   \exp(-2\pi i k\langle\beta',v\rangle-\pi
   ik\tau\langle\beta',\beta'\rangle) 
     \chi^{\sigma_c}_{\lambda,k}(h;\widetilde q)
$$
for any $\beta,\beta'\in \Lambda(T_{w_c})$.
From now on, we shall switch between the different viewpoints for the
characters freely.

\subsection{The action of $\sigma_c$ on the derivation $D$}
Throughout this section, 
let $k$ be a multiple of the basic level of $G$ and let $\lambda\in
P_+^k$ be a highest weight of $\g$ which is invariant under
$\sigma_c$. So the $\sigma_c$-orbit $I$ through $\lambda$ consists of
a single element. 

Let $h_1,\ldots,h_l$ be an orthonormal basis of $\h$ and 
choose $e_\alpha\in\g_\alpha$ and
$f_\alpha\in\g_{-\alpha}$ for each $\alpha\in\Delta_+$ such that the
set $\{e_\alpha,f_\alpha,h_j~|~\alpha\in\Delta_+\text{ and }1\leq
j\leq l\}$ is an orthonormal basis of $\g$. We can make this choice in
such a way 
that $w_c(e_{\alpha})=\pm e_{w_c(\alpha)}$ 
for all
$\alpha\in\Delta$ and similarly for the
$f_{\alpha}$ (see section \ref{adjoint}).

For
$\widetilde\alpha=\alpha+n\delta\in\widetilde\Delta_+^{re}$ let us
define 
$$ 
  e_{\widetilde\alpha} =
  \begin{cases}
  & e_\alpha\otimes z^n\quad\text{ if }
    \alpha\in\Delta_+\,,~n\geq0\,,\qquad\text{ and } \\
  & f_\alpha\otimes z^n\quad\text{ if }
   \alpha\in\Delta_-\,,~n>0\,. 
  \end{cases}
$$
Similarly for
$\widetilde\alpha=\alpha+n\delta\in\widetilde\Delta_+^{re}$, 
we set 
$$
  f_{\widetilde\alpha} =
  \begin{cases}
  &
   f_\alpha\otimes z^{-n}\quad\text{ if }
    \alpha\in\Delta_+\,,~n\geq0\,,
    \qquad\text{ and }\\
  & 
    e_\alpha\otimes z^{-n}\quad\text{ if }
    \alpha\in\Delta_-\,,~n>0.
  \end{cases}
$$    
We use the usual notation $X\otimes z^n=X^{(n)}$ for $X\in\g$.
Then, using the explicit expression
of the Kac-Casimir operator (\cite{K}, 12.8) 
on $V_{\lambda,k}$, we can express $D$ as
an operator on the highest weight representation $V_{\lambda,k}$ of
$\widetilde L(\g)$ via
\begin{equation}\label{D}
  2(k+h^\vee)D = c(\lambda)Id - 2 
   \sum_{\widetilde\alpha\in\widetilde\Delta_+^{re}}
   f_{\widetilde\alpha} e_{\widetilde\alpha} - 
    2\sum_{n=1}^{\infty}\sum_{j=1}^l h_j^{(-n)}h_j^{(n)} - 
      \sum_{j=1}^l (h_j)^2 - \sum_{\alpha\in\Delta_+}h_\alpha\,,
\end{equation}
where $c(\lambda)=\Vert\lambda+\rho\Vert^2-\Vert\rho\Vert^2$ and
$h^\vee$ denotes the dual Coxeter number of $\g$.

In this section, we determine how the expression (\ref{D}) for
the derivation $D$ behaves under the automorphism $\sigma_c$ of
$\widetilde L(\g)$.
Let $(\widetilde\Delta^{re_+})^{\sigma_c}$ be the set of
$\sigma_c$-orbits in $\widetilde\Delta_+^{re}$. Let
$m_{\widetilde\alpha}=|[\widetilde\alpha]|$ denote the cardinality of
the $\sigma_c$-orbit through $\widetilde\alpha$. Let
$\epsilon_{\widetilde\alpha}$ be a
$m_{\widetilde\alpha}$-th root of
unity. Then we set
$$
  e_{\widetilde\alpha}^{\epsilon_{\widetilde\alpha}} = 
    \sum_{j=1}^{m_{\widetilde\alpha}} \epsilon_{\widetilde\alpha}^j
      \sigma_c^j(e_{\widetilde\alpha})\,,
$$
and accordingly 
$$
  f_{\widetilde\alpha}^{\epsilon_{\widetilde\alpha}} = 
    \sum_{j=1}^{m_{\widetilde\alpha}} \epsilon_{\widetilde\alpha}^j
      \sigma_c^j(f_{\widetilde\alpha})\,.
$$
Obviously, $e_{\widetilde\alpha}^{\epsilon_{\widetilde\alpha}}$ and
$f_{\widetilde\alpha}^{\epsilon_{\widetilde\alpha}}$ lie in the
$\epsilon_{\widetilde\alpha}^{-1}$--eigenspace of $\sigma_c$.
Similarly, recall that $\sigma_c$ acts linearly on the space 
$\h\oplus \C C$. For
$h\in\h\oplus\C C$ 
and any $\ord(\sigma_c)$-th root of unity $\epsilon$
we define 
$$h^\epsilon=\sum_{j=1}^{\ord(\sigma_c)}\epsilon^j\sigma_c^j(h)\,.$$ 

Finally, we can decompose $\h=\h_0\oplus\ldots\oplus\h_{\ord(w_c)-1}$
where $\h_j$ denotes the 
$e^{2\pi i \frac{j}{ord(\sigma_c)}}$--eigenspace 
of the action
of the Weyl group element $w_c$ acting on $\h$.
Furthermore, we can choose the basis $h_1,\ldots h_l$ of $\h$ such that
it is adapted to this decomposition of $\h$. 
That is, let $h_{j,1},\ldots h_{j,{l_j}}$ be
an orthonormal basis of the space $\h_j$. Then the set
$\{h_{j,s}~|~0\leq j\leq\ord(w_c)-1, 1\leq s\leq l_j\}$ is an
orthonormal basis of $\h$.

Let us fix a representative $\widetilde\alpha$ for each
$\sigma_c$-orbit $[\widetilde\alpha]$ in $\widetilde\Delta^{re}_+$.
Furthermore, from now on let us fix the roots of unity
$\epsilon_{\widetilde\alpha}=e^{2\pi i\frac{1}{m_{\widetilde\alpha}}}$
and $\epsilon = e^{2\pi i \frac{1}{ord(\sigma_c)}}$
Using the expression from equation (\ref{D}) for the derivation $D$
and the fact that $\sigma_c$ acts as an automorphism
on the universal enveloping algebra of $\widetilde L(\g)$, we can write
\begin{multline}\label{sigD}
  \frac{2(k+h^\vee)}{\ord(\sigma_c)}
    \sum_{j=1}^{\ord(\sigma_c)}\sigma_c^j(D) = \\
      c(\lambda) -
    2 \sum_{[\widetilde\alpha]\in(\widetilde\Delta_+^{re})^{\sigma_c}}
    \frac{1}{m_{\widetilde\alpha}}\sum_{j=1}^{m_{\widetilde\alpha}}
      f_{\widetilde\alpha}^{\epsilon^j_{\widetilde\alpha}}
       e_{\widetilde\alpha}^{\epsilon^{-j}_{\widetilde\alpha}} 
     - \frac{2}{\ord(\sigma_c)} \sum_{n=1}^\infty
    \sum_{j,r=1}^{\ord(w_c)}\sum_{s=1}^{l_j}\epsilon^{2rj}
    {h_{j,s}}^{(-n)}{h_{j,s}}^{(n)} \\
    -\frac{1}{\ord(\sigma_c)^2}\sum_{j=1}^{\ord(\sigma_c)}\sum_{s=1}^{l_j}
    \sum_{r=1}^{\ord(\sigma_c)}h_{j,s}^{\epsilon^r}h_{j,s}^{\epsilon^{-r}}
    - \frac{1}{\ord(\sigma_c)}
    \sum_{\alpha\in\Delta_+}\sum_{j=1}^{\ord(\sigma_c)}\sigma^j(h_\alpha)
\end{multline}


\subsection{The differential equation}\label{theDE}
We now derive a differential equation for the character
$\chi_{\lambda,k}$.

Recall that $\sigma_c$ acts on $\widetilde\h$.
Fix some $h\in\h^{w_c}$ and some $\tau\in\C$. Then the element
$$
  \widetilde H_0 = 
   h-\tau\frac{1}{\ord(\sigma_c)}\sum_{j=1}^{\ord(\sigma_c)}\sigma_c^j(D)
$$
is invariant under $\sigma_c$. But $\widetilde H_0$ contains a non-zero
central term which comes from
$\sum_{j=1}^{\ord(\sigma_c)}\sigma_c^j(D)$. 
Let $C(D)$ denote this central term. Then 
$$
  H_0=\widetilde H_0 -C(D) \in\h^{w_c}\oplus\C D\,,
$$
and $H_0$ is invariant under the residual action of $\sigma_c$ on
$\h\oplus\C D\cong\widetilde\h / \C C$.
Therefore we have
$$
  e^{2\pi i H_0} = e^{2\pi i h}
   q^{(-\frac{1}{ord(\sigma_c)}\sum_{j=1}^{\ord(\sigma_c)}\sigma_c^j(D))} 
    q^{C(D)}
    \in (T\times\C^*)^{\sigma_c}\,
$$
where we have set $q=e^{2\pi i\tau}$ as usual. If $Im(\tau)>0$, we have
$|q|<1$ so that the character 
$\chi_{\lambda,k}(\sigma_ce^{\pi i H_0})$ converges. (Remember that
we have identified the character $\chi_{\lambda,k}$ 
with a section of the family of line bundle $\L^k_{w_c}$ 
which we view locally as a function on $(T\times\C^*)^{\sigma_c}_0$.)

We now want to calculate the derivative
$$\frac{2(k+h^\vee)}{2\pi i}\frac{\partial}{\partial\tau}
\chi_{\lambda,k}(\sigma_c e^{2\pi i H_0})\,.$$  
For simplicity, we will calculate 
$\frac{2(k+h^\vee)}{2\pi i}\frac{\partial}{\partial\tau}
\chi_{\lambda,k}(\sigma_c e^{2\pi i \widetilde H_0})$ and
disregard all central terms which come from the action of
$\sigma_c$ on the derivation $D$.  
Using equation (\ref{sigD}), we can compute
\begin{multline}\label{dtau}
  \frac{2(k+h^\vee)}{2\pi i}\frac{\partial}{\partial\tau}
    \chi_{\lambda,k}(\sigma_c e^{2\pi i \widetilde H_0}) = 
       - 
       Tr_{V_{\lambda,k}}\big(\sigma_ce^{2\pi i \widetilde H_0}
       \frac{2(k+h^\vee)}{\ord(\sigma_c)}
        \sum_{j=1}^{\ord(\sigma_c)}\sigma_c^j(D)\big)\\
    =   -c(\lambda) \chi_{\lambda,k}(\sigma_c e^{2\pi i \widetilde H_0})
    +2\sum_{[\widetilde\alpha]\in(\widetilde\Delta^{re}_+)^{\sigma_c}}
      \frac{1}{m_{\widetilde\alpha}}A_{[\widetilde\alpha]} \\
    +\frac{2}{\ord(\sigma_c)}\sum_{n=1}^\infty\sum_{j=1}^{\ord(w_c)} B_{n,j}\\
     +\frac{1}{\ord(\sigma_c)^2}\sum_{j=1}^{\ord(\sigma_c)} C_{j}
       +\frac{1}{\ord(\sigma_c)}
      \sum_{\alpha\in\Delta_+}\sum_{j=1}^{\ord(\sigma_c)}E_{\alpha,j}\,,
\end{multline}
where we have set 
$$
  A_{[\widetilde\alpha]} = 
    \sum_{j=1}^{m_{\widetilde\alpha}} A_{\widetilde\alpha,j}
$$
with
$$
A_{\widetilde\alpha,j} = Tr_{V_{\lambda,k}}\big(\sigma_ce^{2\pi i
  \widetilde H_0}
    f_{\widetilde\alpha}^{\epsilon^j_{\widetilde\alpha}}
       e_{\widetilde\alpha}^{\epsilon^{-j}_{\widetilde\alpha}} \big)\,
$$
for some fixed $\widetilde\alpha\in[\widetilde\alpha]$. It is clear
that $A_{[\widetilde\alpha]}$ does not depend on the choice of the
representative $\widetilde\alpha$. Similarly,
$$
  B_{n,j}=\sum_{s=1}^{\dim(\h_j)}\sum_{r=1}^{\ord(w_c)} \epsilon^{2jr}
  B_{n,j,s}
$$
with
$$
   B_{n,j,s} = Tr_{V_{\lambda,k}}\big(\sigma_ce^{2\pi i \widetilde H_0}
   h_{s,j}^{(-n)}h_{s,j}^{(n)}\big)\,,
$$
and 
$$
  C_j=\sum_{s=1}^{l_j}\sum_{r=1}^{\ord(\sigma_c)}C_{j,r,s}
$$
with
$$
  C_{j,s,r} = Tr_{V_{\lambda,k}}\big(\sigma_ce^{2\pi i \widetilde H_0}
    h_{j,s}^{\epsilon^r}h_{j,s}^{\epsilon^{-r}}\big)\,.
$$
Finally,
$$
  E_{\alpha,j} = Tr_{V_{\lambda,k}}\big(\sigma_ce^{2\pi i \widetilde H_0}
     \sigma_c^j(h_\alpha)\big)\,.   
$$

We can compute the summands of
equation (\ref{dtau}) more explicitly. Let us first note that since
$\widetilde H_0$ is invariant under $\sigma_c$, we have
$\sigma_c^j(\widetilde\alpha)(\widetilde
H_0)=\widetilde\alpha(\widetilde H_0)$ for all
$j\in\N$, so that for any root $\widetilde\alpha\in\widetilde\Delta$, 
the value of $\widetilde H_0$ on the
$\sigma_c$-orbit through 
$\widetilde\alpha$ is constant.
Furthermore, for any real root $\widetilde\alpha\in\widetilde\Delta^{re}$,
we have chosen the $e_{\widetilde\alpha}$ and
$f_{\widetilde\alpha}$ such that 
$\sigma_c(e_{\widetilde\alpha})=\pm e_{\widetilde\alpha}$ 
and similarly for $f_{\widetilde\alpha}$. Let us define
\begin{equation}\label{salpha}
  s(\widetilde\alpha)=
\begin{cases} 
  -1&
  \text{if} \quad\sigma_c(\widetilde\alpha)= \widetilde\alpha
  \text{ and }\sigma_c(e_{\widetilde\alpha})=-e_{\widetilde\alpha} \\
  1&
  \text{in all other cases}
\end{cases}
\end{equation}
If $s(\widetilde\alpha)=1$, 
we can calculate
\begin{align*}
  A_{\widetilde\alpha,j}&  =  Tr_{V_{\lambda,k}}
   \big(\sigma_ce^{2\pi i \widetilde H_0}
    f_{\widetilde\alpha}^{\epsilon^j_{\widetilde\alpha}}
     e_{\widetilde\alpha}^{\epsilon^{-j}_{\widetilde\alpha}} \big) \\
  & = 
    Tr_{V_{\lambda,k}}\big(
  e_{\widetilde\alpha}^{\epsilon^{-j}_{\widetilde\alpha}}\sigma_c
   e^{2\pi i \widetilde H_0}
    f_{\widetilde\alpha}^{\epsilon^j_{\widetilde\alpha}}\big)\\
  & = 
    \epsilon_{\widetilde\alpha}^{-j} 
   Tr_{V_{\lambda,k}}\big(
   \sigma_c e_{\widetilde\alpha}^{\epsilon^{-j}_{\widetilde\alpha}}
   e^{2\pi i \widetilde H_0}
    f_{\widetilde\alpha}^{\epsilon^j_{\widetilde\alpha}}\big)\\
  & = 
    \epsilon_{\widetilde\alpha}^{-j} 
    e^{-2\pi i \widetilde\alpha(\widetilde H_0)}
     Tr_{V_{\lambda,k}}\big(
      \sigma_c e^{2\pi i \widetilde H_0} 
       e_{\widetilde\alpha}^{\epsilon^{-j}_{\widetilde\alpha}}
        f_{\widetilde\alpha}^{\epsilon^j_{\widetilde\alpha}}\big)\\
      & = 
    \epsilon_{\widetilde\alpha}^{-j} 
    e^{-2\pi i \widetilde\alpha(\widetilde H_0)}
     Tr_{V_{\lambda,k}}\big(
      \sigma_c e^{2\pi i \widetilde H_0} 
     (f_{\widetilde\alpha}^{\epsilon^j_{\widetilde\alpha}}
     e_{\widetilde\alpha}^{\epsilon^{-j}_{\widetilde\alpha}} + 
     [e_{\widetilde\alpha}^{\epsilon^{-j}_{\widetilde\alpha}},
       f_{\widetilde\alpha}^{\epsilon^j_{\widetilde\alpha}}])\big).
\end{align*}
In the first step, we have used the cyclic property of the trace. In
the second step we have used the fact, that
$e_{\widetilde\alpha}^{\epsilon^{-j}_{\widetilde\alpha}}$ lies in the
$\epsilon_{\widetilde\alpha}^j$--eigenspace of the action of $\sigma_c$ on 
$\widetilde L(\g)$. 
The remaining steps are the standard commutation relations in
$\widetilde L(\g)$. The commutator in the last step can be computed as
follows:
\begin{align*}
  [e_{\widetilde\alpha}^{\epsilon^{-j}_{\widetilde\alpha}},
       f_{\widetilde\alpha}^{\epsilon^j_{\widetilde\alpha}}] = & 
   \sum_{s=1}^{m_{\widetilde\alpha}}\sum_{r=1}^{m_{\widetilde\alpha}}
       \epsilon^{-js}_{\widetilde\alpha}\epsilon^{jr}_{\widetilde\alpha}
       [\sigma_c^s(e_{\widetilde\alpha}),\sigma_c^r(f_{\widetilde\alpha})] \\
    = & 
     \sum_{s=1}^{m_{\widetilde\alpha}}
      \sigma_c^s([e_{\widetilde\alpha},f_{\widetilde\alpha}])\,,
\end{align*}
since we have
$[\sigma_c^s(e_{\widetilde\alpha}),\sigma_c^r(f_{\widetilde\alpha})]=0$
whenever $s\neq r$. Indeed, by assumption, we have 
$\sigma_c^s(\widetilde\alpha)\neq\sigma_c^r(\widetilde\alpha)$ for
$s\neq r$. Since $\sigma_c$ acts as a diagram automorphism on the
Dynkin diagram of the affine root system $\widetilde\Delta$, there
exists a root basis $\widetilde\Pi'$ 
of $\widetilde\Delta$ which contains both  
$\sigma_c^s(\widetilde\alpha)$ and
$\sigma_c^r(\widetilde\alpha)$. Now, we can find some element $w$ in
the Weyl group of $\widetilde\Delta$ such that
$w(\widetilde\Pi')=\widetilde\Pi$ (\cite{K}, Proposition 5.9). So we
have $[\sigma_c^s(e_{\widetilde\alpha}),\sigma_c^r(f_{\widetilde\alpha})]=
[e_{\widetilde\alpha_i},f_{\widetilde\alpha_j}]$ for two
simple roots
$\widetilde\alpha_i\neq\widetilde\alpha_j\in\widetilde\Pi$. 
The last commutator has to vanish due to the Serre relations for
$\widetilde L(\g)_\pol$ (see \cite{K}, 1.2).

If $\widetilde\alpha=\alpha+n\delta$ with $\alpha\in\Delta_+$, we have 
$$
  [e_{\widetilde\alpha},f_{\widetilde\alpha}]=h_\alpha+nC\,,
$$
where $h_\alpha\in\h$ is the co-root of $\g$ corresponding to the root
$\alpha\in\Delta$.
Similarly, if  $\widetilde\alpha=\alpha+n\delta$ with $\alpha\in\Delta_-$, 
one has 
$$
  [e_{\widetilde\alpha},f_{\widetilde\alpha}]=-h_\alpha+nC\,.
$$
Let us write $n(\widetilde\alpha)=n$ for
$\widetilde\alpha=\alpha+n\delta$ and $h_{\widetilde\alpha}=h_\alpha$
if $\alpha\in\Delta_+$ and $h_{\widetilde\alpha}=-h_\alpha$ if 
$\alpha\in\Delta_-$.
Then the calculations above give us
\begin{multline*}
    A_{\widetilde\alpha,j} = 
     \epsilon_{\widetilde\alpha}^{-j} 
     e^{-2\pi i \widetilde\alpha(\widetilde H_0)}
       \big(A_{\widetilde\alpha,j} ~+ \\
            Tr_{V_{\lambda,k}}\big(
              \sigma_c e^{2\pi i \widetilde H_0} 
      \sum_{s=1}^{m_{\widetilde\alpha}}
           \sigma_c^s(h_{\widetilde\alpha}) + 
        \chi_{\lambda,k}(\sigma_c e^{2\pi i \widetilde H_0})\big)
          k m_{\widetilde\alpha} n(\widetilde\alpha)\big)\,,
\end{multline*}
or equivalently
\begin{multline*}
  A_{\widetilde\alpha,j} = 
  \frac{\epsilon_{\widetilde\alpha}^{-j} 
    e^{-2\pi i \widetilde\alpha(\widetilde H_0)}}
      {1-\epsilon_{\widetilde\alpha}^{-j} 
    e^{-2\pi i \widetilde\alpha(\widetilde H_0)}}
     \big(
       Tr_{V_{\lambda,k}}\big(
              \sigma_c e^{2\pi i \widetilde H_0} 
      \sum_{s=1}^{m_{\widetilde\alpha}}
           \sigma_c^s(h_{\widetilde\alpha})\big) \\
   + 
    \chi_{\lambda,k}(\sigma_c e^{2\pi i \widetilde H_0})
         k m_{\widetilde\alpha} n(\widetilde\alpha)\big)
\end{multline*}
Summing over the $j$'s, we can write
\begin{multline*}
  A_{[\widetilde\alpha]} = 
     \frac{
      m_{\widetilde\alpha}
      e^{-2\pi i m_{\widetilde\alpha}\widetilde\alpha(\widetilde H_0)}}
      {1-e^{-2\pi i m_{\widetilde\alpha}\widetilde\alpha(\widetilde H_0)}}
      \big(Tr_{V_{\lambda,k}}\big(
              \sigma_c e^{2\pi i \widetilde H_0} 
     \sum_{s=1}^{m_{\widetilde\alpha}}
           \sigma_c^s(h_{\widetilde\alpha})\big)\\
   + 
    \frac{k m_{\widetilde\alpha}^2 n(\widetilde\alpha)
     e^{-2\pi i m_{\widetilde\alpha}\widetilde\alpha(\widetilde H_0)}}
    {1-e^{-2\pi i m_{\widetilde\alpha}\widetilde\alpha(\widetilde H_0)}}
    \chi_{\lambda,k}(\sigma_c e^{2\pi i \widetilde H_0})\big)\,.
\end{multline*}
Since 
$$
  \sum_{s=1}^{m_{\widetilde\alpha}}
           \sigma_c^s(h_{\widetilde\alpha}) = \sum_{s=1}^{m_{\widetilde\alpha}}
           w_c^s(h_{\widetilde\alpha}) + y_{\widetilde\alpha} C
$$
for some $y_\alpha\in\C$,
we can finally write
\begin{multline}\label{A}
  A_{[\widetilde\alpha]} = 
     \frac{
        m_{\widetilde\alpha}
      e^{-2\pi i m_{\widetilde\alpha}\widetilde\alpha(\widetilde H_0)}}
      {1-e^{-2\pi i m_{\widetilde\alpha}\widetilde\alpha(\widetilde H_0)}}
           \sum_{s=1}^{m_{\widetilde\alpha}} 
         \frac{1}{2\pi i}\frac{\partial}{\partial w_c^s(h_{\widetilde\alpha})}
           \chi_{\lambda,k}(\sigma_c e^{2\pi i \widetilde H_0})
   + \\
    \frac{k m_{\widetilde\alpha}^2 n(\widetilde\alpha)
     e^{-2\pi i m_{\widetilde\alpha}\widetilde\alpha(\widetilde H_0)}}
    {1-e^{-2\pi i m_{\widetilde\alpha}\widetilde\alpha(\widetilde H_0)}}
    \chi_{\lambda,k}(\sigma_c e^{2\pi i \widetilde H_0})
    + Y_{[\widetilde\alpha]}C
\end{multline}
for some $Y_{[\widetilde\alpha]}\in\C$ 
which comes from the action of $\sigma_c$ on the
derivation $D$.

If $s(\widetilde\alpha)=-1$, we necessarily have
$m_{\widetilde\alpha}=1$. A similar calculation as above gives
\begin{multline}\label{As}
  A_{[\widetilde\alpha]} = 
     \frac{ 
      -e^{-2\pi i \widetilde\alpha(\widetilde H_0)}}
      {1+e^{-2\pi i \widetilde\alpha(\widetilde H_0)}} 
         \frac{1}{2\pi i}
           \frac{\partial}{\partial w_c^s(h_{\widetilde\alpha})}
           \chi_{\lambda,k}(\sigma_c e^{2\pi i \widetilde H_0})
   - \\
    \frac{k n(\widetilde\alpha)
     e^{-2\pi i \widetilde\alpha(\widetilde H_0)}}
    {1+e^{-2\pi i \widetilde\alpha(\widetilde H_0)}}
    \chi_{\lambda,k}(\sigma_c e^{2\pi i \widetilde H_0})
    + Y_{[\widetilde\alpha]}C\,.
\end{multline}

A similar calculation for the $B_{n,j,s}$ gives
\begin{align*}
  B_{n,j,s} = & Tr_{V_{\lambda,k}}\big(\sigma_c e^{2\pi i \widetilde H_0}
  {h_{j,s}}^{(-n)}{h_{j,s}}^{(n)} \big) \\
    = &
   \epsilon^{-j}e^{-2\pi i n\delta(\widetilde H_0)}
  Tr_{V_{\lambda,k}}\big(\sigma_c e^{2\pi i \widetilde H_0}
  ({h_{j,s}}^{(-n)}{h_{j,s}}^{(n)}+
  [{h_{j,s}}^{(n)},{h_{j,s}}^{(-n)}])\big)\,. 
\end{align*}
Using the fact that $\delta(\widetilde H_0)=-\tau$ and 
$[{h_{j,s}}^{(n)},{h_{j,s}}^{(-n)}] = n C$, we get
\begin{equation}\label{Bnjs}
    B_{n,j,s} = \frac{\epsilon^{-j}q^n n k
    \chi_{\lambda,k}(\sigma_ce^{2\pi i \widetilde H_0})}
    {1-\epsilon^{-j}q^n}
\end{equation}
This shows that $B_{n,j,s}$ in fact does not depend on $s$. So we get
\begin{equation}
  B_{n,j} =  
    \dim(\h_j)\sum_{r=1}^{\ord(\sigma_c)}\epsilon^{2jr}
    \frac{\epsilon^{-j}q^n n k 
    \chi_{\lambda,k}(\sigma_ce^{2\pi i \widetilde H_0})}
    {1-\epsilon^{-j}q^n} \,.
\end{equation}
This sum vanishes if $\ord(\sigma_c)\not|2j$, otherwise we get
\begin{equation}
  B_{n,j} = 
    \ord(\sigma_c)\dim(\h_j)
    \frac{\epsilon^{-j}q^n n k 
    \chi_{\lambda,k}(\sigma_ce^{2\pi i \widetilde H_0})}
    {1-\epsilon^{-j}q^n} \,.
\end{equation}    

Finally, since $h_{j,s}^{\epsilon^r}$ is in the 
$\epsilon^{-r}$--eigenspace of the
action of $\sigma_c$ on $\h\oplus \C C$, a similar calculation gives
$$
  C_{j,s,r}= Tr_{V_{\lambda,k}}\big(\sigma_ce^{2\pi i \widetilde H_0}
    h_{j,s}^{\epsilon^r}h_{j,s}^{\epsilon^{-r}}\big) = 
    \epsilon^{-r}C_{j,s,r}\,.
$$
This shows that $C_{j,s,r}=0$ unless $\epsilon^r=1$. Now, if $j=0$, we
have $h_{j,s}^{1}=\ord(\sigma_c)h_{j,s}$, so that we get
\begin{equation}
  C_0 = 
     \sum_{s=1}^{l_0}C_{0,s,0} = 
    \frac{\ord(\sigma_c)^2}{(2\pi i)^2}
   \Delta_{\h_0}\chi_{\lambda,k}(\sigma_ce^{2\pi i \widetilde H_0})\,,
\end{equation}
where $\Delta_{h_0}$ denotes the Laplace operator on $\h_0$.
On the other hand, for $j\neq0$, we have
$h_{j,s}^{1}=y_{s,j} C$ for some $y_{s,j}\in\C$ so that we can drop
these terms since they come from the action of $\sigma_c$ on $D$.

Since all roots $\widetilde\alpha$ vanish on the center $C$ of
$\widetilde L(\g)$, we have 
$\widetilde\alpha(\widetilde H_0)=\widetilde\alpha(H_0)$ for all
$\widetilde\alpha\in\widetilde\Delta$. 
So putting the calculations above together, reordering the terms
slightly, and keeping in mind that we have to disregard all central terms
coming from the action of $\sigma_c$ on the derivation $D$, 
we get
\begin{multline}\label{first}
   \frac{2(k+h^\vee)}{2\pi i}\frac{\partial}{\partial\tau}\big(
    \chi_{\lambda,k}(\sigma_c e^{2\pi i H_0})\big) = \\
  \big(-c(\lambda) + A_1+A_2 + B + E  -  
    \frac{1}{4\pi^2}\Delta_{\h_0}\big)
     \chi_{\lambda,k}(\sigma_c e^{2\pi i H_0})
\end{multline}
where we have set 
\begin{align*}
  A_1= & \sum_{[\widetilde\alpha]\in(\widetilde\Delta^{re}_+)^{\sigma_c}}
   \frac{s(\widetilde\alpha)2 
      e^{-2\pi i m_{\widetilde\alpha}\widetilde\alpha(H_0)}}
      {1-s(\widetilde\alpha)
       e^{-2\pi i m_{\widetilde\alpha}\widetilde\alpha(H_0)}}
           \sum_{s=1}^{m_{\widetilde\alpha}} 
    \frac{1}{2\pi i}\frac{\partial}{\partial w_c^s(h_{\widetilde\alpha})}\,,\\
  A_2 = & \sum_{[\widetilde\alpha]\in(\widetilde\Delta^{re}_+)^{\sigma_c}}
    \frac{s(\widetilde\alpha)2 k  m_{\widetilde\alpha}n(\widetilde\alpha)
     e^{-2\pi i m_{\widetilde\alpha}\widetilde\alpha(H_0)}}
    {1-s(\widetilde\alpha)
     e^{-2\pi i m_{\widetilde\alpha}\widetilde\alpha(H_0)}}\,, \\
  B= &\sum_{n=1}^\infty \sum_{j=1}^{\ord(\sigma_c)}
   2 \dim(\h_j)
    \frac{\epsilon^{-j}q^n n k} 
     {1-\epsilon^{-j}q^n}\,, \\
  E = & 
  \frac{1}{\ord(\sigma_c)}\sum_{\alpha\in\Delta_+}\sum_{j=1}^{\ord(\sigma_c)}
  \frac{1}{2\pi i}\frac{\partial}{\partial w_c^j(h_\alpha)}\,.
\end{align*}

Equation (\ref{first}) is the desired differential equation for the
characters. We will now simplify the equation. To do this, let us
make some definitions. First, recall the element
$\rho=\frac12\sum_{\alpha\in\Delta_+}\alpha$ and set
$$
  \rho_{w_c}=\frac{1}{\ord(w_c)}\sum_{j=1}^{\ord(w_c)}w_c(\rho)\,.
$$ 
Then we can define the function
$$
  \F_{w_c}(H_0)=
    e^{2\pi i \rho_{w_c}(H_0)}
    \prod_{[\widetilde\alpha]\in(\widetilde\Delta^{re}_+)^{\sigma_c}}
   \big(1-s(\widetilde\alpha)
    e^{-2\pi i m_{\widetilde\alpha}\widetilde\alpha(H_0)}\big)
    \prod_{n=1}^\infty
    \prod_{j=1}^{\ord(w_c)}(1-\epsilon^{-j}q^n)^{\dim(\h_j)}\,.
$$
By the definition of $H_0$ and invariance of $\delta$ under
$\sigma_c$, we have $\widetilde\alpha(H_0) = \alpha(H_0) - \tau n$ 
for $\widetilde\alpha=\alpha+n\delta$.
Hence we get
$$
  \frac{1}{2\pi i}\frac{\partial}{\partial\tau} \F_{w_c}(H_0)
    = - \big(\frac{1}{2k} A_2
    + \frac{1}{2k}B\big)\F_{w_c}(H_0)\,.
$$
Now, $\F_{w_c}$ is a section of the line bundle $\L_{w_c}^{h^\vee}$
from remark \ref{Lk2}, so
$\chi_{\lambda,k}\F_{w_c}$ is a section of
$\L_{w_c}^{k+h^\vee}$. Therefore we get
\begin{multline}
  \frac{1}{\F_{w_c}}
  \frac{2(k+h^\vee)}{2\pi i}\frac{\partial}{\partial\tau}
   \big(\chi_{\lambda,k}\F_{w_c}\big) = \\
     \frac{2(k+h^\vee)}{2\pi i}
    \left(\frac{\partial \chi_{\lambda,k}}{\partial\tau} + 
  \frac{\chi_{\lambda,k}}{\F_{w_c}} 
    \frac{\partial \F_{w_c}}{\partial\tau} \right)\\
   = 
    \big(-c(\lambda) 
    + A_1 + E - \frac{1}{4\pi^2}\Delta_{h^{w_c}}\big)\chi_{\lambda,k}
\end{multline}

Next, we take care of the Laplacian. We use the formula  
$$ 
  \Delta_{\h^{w_c}}(\chi_{\lambda,k}\F_{w_c}) = 
    (\Delta_{\h^{w_c}}\chi_{\lambda,k})\F_{w_c} + 
      2\langle\nabla\chi_{\lambda,k},
       \nabla\F_{w_c}\rangle +
    \chi_{\lambda,k}(\Delta_{\h^{w_c}}\F_{w_c})\,.
$$
But we have
\begin{align}
    \frac{1}{2\pi i}\nabla(\F_{w_c}) = &
      \sum_{s=1}^{\dim(\h_0)}\rho_{w_c}(h_{0,s}) 
      \F_{w_c}h_{0,s} + \\ 
     & 
       \sum_{s=1}^{\dim(\h_0)} 
     \sum_{[\widetilde\alpha]\in(\widetilde\Delta^{re}_+)^{\sigma_c}}
     m_{\widetilde\alpha}\widetilde\alpha(h_{0,s})
    \frac{s(\widetilde\alpha)e^{-2\pi i m_{\widetilde\alpha}
      \widetilde\alpha(H_0)}}
      {(1-s(\widetilde\alpha)
     e^{-2\pi i m_{\widetilde\alpha}\widetilde\alpha(H_0)})} 
     \F_{w_c}h_{0,s} 
     \\ 
     = &
     \F_{w_c} h_{\rho_{w_c}} + 
     \sum_{[\widetilde\alpha]\in(\widetilde\Delta^{re}_+)^{\sigma_c}}
      m_{\widetilde\alpha}
    \frac{s(\widetilde\alpha)e^{-2\pi i m_{\widetilde\alpha}
        \widetilde\alpha(H_0)}}
      {(1-s(\widetilde\alpha)
      e^{-2\pi i m_{\widetilde\alpha}\widetilde\alpha(H_0)})} 
     \F_{w_c}h_{\alpha}
\end{align}
where we have used the notation from before for an orthogonal basis
$h_{0,s}$ of $\h_0=\h^{w_c}$. Since $h_{0,s}\in\h$, we have 
$\widetilde\alpha(h_{0,s})=\alpha(h_{0,s})$ for
$\widetilde\alpha=\alpha+n\delta$. 
Now, since we have
$$
  \langle\nabla\chi_{\lambda,k},h\rangle = 
   \frac{\partial}{\partial h}\chi_{\lambda,k}\,
$$
and 
$$
  h_{\rho_{w_c}} = \frac{1}{2\ord(\sigma_c)}
    \sum_{\alpha\in\Delta_+}\sum_{j=1}^{\ord(w_c)}w_c^j(h_\alpha)\,,
$$
we get
\begin{equation}
  \langle\frac{1}{2\pi i}\nabla\chi_{\lambda,k},
  \frac{1}{2\pi i}\nabla\F_{w_c}\rangle =  
   \frac12 (E\chi_{\lambda,k})\F_{w_c} +
    \frac12 (A_1\chi_{\lambda,k})\F_{w_c}
\end{equation}
So putting everything together, we get
$$
  \frac{1}{\F_{w_c}}\left( \frac{2(k+h^\vee)}{2\pi i}
   \frac{\partial}{\partial\tau}
     \big(\chi_{\lambda,k}\F_{w_c}\big)
      +\frac{1}{4\pi^2}\Delta_{\h_0}(\chi_{\lambda,k}\F_{w_c}\big)
       - \chi_{\lambda,k}\frac{1}{4\pi^2}\Delta_{\h_0}\F_{w_c}\right) 
   = -c(\lambda)\chi_{\lambda,k}\,.
$$
A direct calculation gives 
$$
  \frac{1}{4\pi^2}\Delta_{\h_0}\F_{w_c} = - \vert\rho_{w_c}\vert^2 \F_{w_c}
$$
so that we arrive at the final equation:
\begin{Theorem}\label{thmdgl}
  Let $k$ be a multiple of the basic level of $G$, and let
  $\lambda\in P_+^k$ be a highest weight of $\widetilde G$ which is
  invariant under the action of $\sigma_c$. Then the character 
  $\chi_{\lambda,k}$ of the representation $V_{\lambda,k}$ restricted
  to the torus $(\sigma_c,1)\big(T\times\C^*\big)^{\sigma_c}_0$ has
  to satisfy the following differential equation:
  \begin{equation}
    \frac{1}{\F_{w_c}}\left( - \frac{2(k+h^\vee)}{2\pi
  i}\frac{\partial}{\partial \tau}
    +\frac{1}{(2\pi i)^2}\Delta_{\h_0} - \vert\rho_{w_c}\vert^2 \right)
    \chi_{\lambda,k}\F_{w_c} = c(\lambda)\chi?{lambda,k}
   \end{equation}
\end{Theorem}
 
\begin{Remark}
  If we set  $\sigma_c=id$ theorem \ref{thmdgl} gives a differential
  equation for the characters of highest weight representations of
  loop groups based on simply connected Lie groups. This case 
  is treated in more generality in \cite{EK} (see also \cite{EFK}), 
  where a differential
  equation for traces of 
  intertwining maps between representations of loop groups is
  derived. Introducing the automorphism $\sigma_c$, one can find a
  differential equation for traces of 
  intertwining maps between representations
  of loop groups based on non-simply connected Lie groups. 
\end{Remark}
\begin{Remark}
  We have only treated the case of automorphisms $\sigma_c$
  associated to the center of the simply connected Lie group
  $\widetilde G$. The same calculations work if one replaces the
  automorphism $\sigma_c$ by an arbitrary finite order automorphism $\sigma$ 
  of the Lie
  algebra $\widetilde L(\g)$ which comes from an automorphism of the
  Dynkin diagram of $\widetilde L(\g)$.
\end{Remark}



\section{Character formulas}\label{char}
\subsection{Theta functions}\label{thetafns}
The line bundle $\mathcal{L}^k$ introduced in Remark \ref{Lk}
has been studied in \cite{Lo}. Here we give a brief account of the main
results.
As in the end of section \ref{charsec} fix some 
$\tau\in\C$ with $Im(\tau)>0$ and let
$q=e^{2\pi i\tau}$. Let $L\subset\C$ be the lattice generated by $1$ and 
$\tau$. Then the elliptic curve $E_\tau=\C/L$
is isomorphic to the curve $E_q$ considered in the last paragraph via the 
map $x\mapsto e^{2\pi i x}$. 
Denote by $\Gamma(\L^k)$ the
space of holomorphic sections  
of $\L^k$. Let 
$\Gamma(\L^k)^W$, resp. $\Gamma(\L^k)^{-W}$ denote the spaces of holomorphic
$W$-invariant, resp. $W$-anti invariant sections of $\L^k$. Looijenga's theorem
describes  the set $\Gamma(\L^k)^{-W}$ explicitly in
terms of theta functions. Since we also want to consider sections in the line
bundles $\L^k_{w_c}$ from Remark \ref{Lk2}, we start with some general
statements.

Let $V$ be a finite dimensional Euclidean vector space and denote the
bilinear 
form in $V$ by $\langle~,~\rangle$. In what follows, we will freely identify
$V$ and $V^*$ via the bilinear form $\langle~,~\rangle$. Let $M$ be some
integer lattice in $V$. Let us denote by $M^*\subset V$ 
the dual lattice corresponding to  $M$. That is, $M^*=\{\lambda\in
V~|~\langle\alpha,\lambda\rangle\in\Z\text{ for all }\alpha\in M\}$.
For $\mu\in M^*$ and $k\in\N$, 
define the theta function $\Theta_{\mu,k}$ on $V\otimes\C$ 
via
$$
  \Theta_{\mu,k}(v)=\exp(-\frac{1}{k}\pi i \tau  \langle \mu,\mu\rangle)
 \sum_{\gamma\in \frac{1}{k}\mu+M}
 \exp(2\pi i k\tau(\langle\gamma,\tau^{-1}v\rangle+
 \frac12\langle\gamma,\gamma\rangle))\,.
$$
This function converges absolutely on compact sets and satisfies the identity
\begin{equation}\label{id}
  \Theta_{\mu,k}(v+\beta+\tau\beta')=
   \exp(-2\pi i k\langle\beta',v\rangle-\pi
   ik\tau\langle\beta',\beta'\rangle) 
     \Theta_{\mu,k}(v)
\end{equation}
for any $\beta,\beta'\in \Lambda(T)$.
To emphasize the dependence
of $\Theta_{\mu,k}$ on $\tau$,
we will write $\Theta_{\mu,k}(v;\tau)$.

\medskip

Now suppose that $G$ is a simply connected Lie group and 
take $V_\R=\h_\R$ endowed with the normalized Killing form 
$\langle~,~\rangle$.
Since $\Lambda(T)$ is an integer lattice with respect to $\langle~,~\rangle$, 
we can take $M$ to be
$\Lambda(T)$. It is clear from equation (\ref{id}) 
that in this situation, the functions
$\Theta_{\mu,k}(\,\cdot\,;\tau)$ define holomorphic sections of the line bundle
$\L^k$.

We can define the anti-invariant theta functions
$$
  \sA\Theta_{\mu,k} = \sum_{w\in W}(-1)^{l(w)}\Theta_{w\mu},
$$ 
where $W$ the Weyl group of $G$, and $l(w)$ denotes the length of $w$
in $W$.
\par
Let $\mathfrak{a}\subset\h_\R$ be a fundamental alcove for 
the action of the affine Weyl group $\widetilde W=W\ltimes\Lambda(T)$. 
Then the following result is due to \cite{Lo}. 
\begin{Prop}
\label{thetabasis}
  Let $k$ be a positive integer.
  The anti-invariant theta functions $\sA\Theta_{\mu,k}$ with 
  $\mu\in \Lambda(T)^*\cap k\mathfrak{a}$ 
  form a basis in $\Gamma(\L^{k})^{-W}$.
\end{Prop}

Let $h^\vee$ denote the dual Coxeter number of the root system
$\Delta$ of $\g$.
As usual, we set $\rho=\frac12\sum_{\alpha\in \Delta_+}\alpha$, where
$\Delta_+$
denotes the set of positive roots of $\Delta$ with respect to the Weyl 
chamber containing the fundamental alcove $\mathfrak{a}$.
For $\sigma_c=id$ and $h\in\h$, 
the function $\F=\F_{id}$ which was defined in
section \ref{theDE} can be written as
\begin{multline*}
  \F_{id}(h-\tau D) = \exp(2\pi i \langle\rho,h\rangle)
   \prod_{n=1}^\infty\left(1-q^{n})\right)^l
   \prod_{\alpha\in \Delta_+} 
    \left(1-\exp(-2\pi i \langle\alpha,h\rangle)\right)\\
    \prod_{\alpha\in \Delta}\prod_{n=1}^\infty\left(1-q^n\exp(-2\pi i
     \langle\alpha,z\rangle)\right)\,,
\end{multline*}
where, as before, $l$ denotes the rank of $\g$. It is well known that
$\F$  is a 
holomorphic, $W$ anti-invariant section of $\L^{h^\vee}$ (see \cite{Lo}). 

\medskip

If $G$ is of the form $G=\widetilde G/Z$, where $Z=\langle c\rangle$ 
is a cyclic
subgroup of the center of the simply connected group $\widetilde G$, 
we have to consider
the lattice $\Lambda(T_{w_c})\subset\h^{w_c}$. Let $\langle~,~\rangle$ denote
the normalized Killing form restricted to $\h^{w_c}_\R$. As before, let $k_b$
denote the basic level of $G$. A case by case check shows that
$\Lambda(T_{w_c})$ is an integral lattice with respect to the bilinear form
$k_b\langle~,~\rangle$. 
Denote by
$\Lambda(T_{w_c})^*_{k_b}\subset\h^{w_c}_\R$ the 
dual lattice of 
$\Lambda(T_{w_c})$ with respect to the bilinear form $k_b\langle~,~\rangle$. 
Then, for $\mu\in \Lambda(T_{w_c})^*_{k_b}$ and $k$ a multiple of $k_b$,
the theta function $\Theta_{\mu,k}$ is a section of the line bundle
$\L_{w_c}^k$ which was defined in Remark \ref{Lk2}.
As before, we define the anti-invariant theta functions
$$
  \sA_0\Theta_{\mu,k} = \sum_{w\in W_0}(-1)^{l(w)}\Theta_{w\mu},
$$ 
where $W_0$ is the finite Weyl group introduced in Proposition
\ref{tildeW}, 
and $l(w)$ denotes the length of $w$ in $W_0$.  
Let $\mathfrak{a}_{w_c}\subset\h^{w_c}_\R$ be a fundamental alcove for 
the action of the affine Weyl group 
$\widetilde W^{\sigma_c}=W_0\ltimes\Lambda(T_{w_c})$. 
Then we have (see \cite{K}, Chapter 13):
\begin{Prop}
\label{thetabasis2}
  Let $k$ be a multiple of the basic level $k_b$ of $G$.
  The anti-invariant theta functions $\sA_0\Theta_{\mu,k}$ with 
  $\mu\in \Lambda(T_{w_c})^*_{k_b}\cap k\mathfrak{a}_{w_c}$ 
  form a basis in $\Gamma(\L_{w_c}^{k})^{-W_0}$.
\end{Prop}

\medskip

Finally, the following proposition is straight forward to check.
\begin{Prop}
  The function $\F_{w_c}$ defines a $W_0$--anti-invariant
  holomorphic section of the line bundle $\L^{h^\vee}_{w_c}$.
\end{Prop}


\subsection{The Kac-Weyl character formula}\label{sec:KW}
Throughout this section let $G$ be simply connected.
Let $\chi_{\lambda,k}$ denote the character of the $\widetilde L(G)$-module
$V_{\lambda,k}$ of highest weight $\lambda$ and level $k$. The goal of
this section is to find an explicit formula for the character 
$\chi_{\lambda,k}$ viewed as a section of the family of line bundles
$\L^k$ as described in Section \ref{thetafns}.

Since $\F=\F_{id}$ is a $W$-anti-invariant section of $\L^{h^\vee}$,
the product $\chi_{\lambda,k}\F$
defines a $W$-anti-invariant section of the bundle $\L^{(k+h^\vee)}$.
By Proposition (\ref{thetabasis}), we can write this product uniquely as
\begin{equation}\label{uniquely}
  \chi_{\lambda,k}(h;\tau)\F(h;\tau) = 
   \sum_{\mu\in\Lambda(T)^*\cap (k+h^\vee)\mathfrak{a}}
   f_\mu(\tau)\sA\Theta_{\mu,k+h^\vee}(h;\tau)\,. 
\end{equation}
We can naturally identify the lattice  $\Lambda(T)^*$, with the weight
lattice $P$ of $\g$. So the sum  in equation (\ref{uniquely}) ranges
over $P\cap (k+h^\vee)\mathfrak{a}$.

Now we let $\tau$ vary in the upper half plane. In the situation at
hand, the differential equation from theorem \ref{thmdgl} reads
\begin{multline}\label{DGL}
  \left(-\frac{2(k+h^\vee)}{2\pi i}\frac{\partial}{\partial\tau} +
  \frac{1}{(2\pi i)^2}\Delta_\h 
  - \langle\rho,\rho\rangle\right)(\chi_{\lambda,k}(h,\tau)\F(h,\tau)) = \\
  \langle\lambda,\lambda + 2\rho\rangle
  \chi_{\lambda,k}(h,\tau)\F(h,\tau),
\end{multline}
where $\Delta_\h$ is the Laplacian on $\h$.

Substituting equation (\ref{uniquely}) into differential equation (\ref{DGL}) 
and keeping in mind that the 
$\sA\Theta_{\mu,k}$ are linearly independent, we find 
\begin{equation}\label{fmu} 
  f_\mu(\tau)= a_\mu\exp\left(\frac{2\pi i \tau}{2(k+h^\vee)}
  (\langle\lambda+\rho,\lambda+\rho\rangle - 
  \langle\mu,\mu\rangle)\right)\,,
\end{equation}
where $a_\mu$ is a constant depending only on $\mu$. So we get
\begin{multline}\label{sum}
  \chi_{\lambda,k}(h,\tau)\F(h,\tau)=\\
  \sum_{\mu\in P\cap (k+h^\vee)\mathfrak{a}} 
  a_\mu\exp\left(\frac{2\pi i \tau}{2(k+h^\vee)}
  (\langle\lambda+\rho,\lambda+\rho\rangle - 
  \langle\mu,\mu\rangle)\right)
  \sA\Theta_{\mu,k+h^\vee}(h;\tau).
\end{multline}

\medskip

By definition of $\chi_{\lambda,k}$, we can write
\begin{equation}\label{ch}  
  \chi_{\lambda,k}(h,\tau) =
  \sum_{\tilde\mu\in P(\lambda,k)}\dim V_{\lambda,k}[\tilde\mu]\cdot 
  q^{-D(\tilde\mu)} \exp(2\pi i\langle\tilde\mu,h\rangle),
\end{equation}
where $P(\lambda,k)\subset(\h\oplus\C C\oplus\C D)^*$ 
denotes the set of weights of $V_{\lambda,k}$, and for
$\tilde\mu\in P(\lambda,k)$, the space 
$V_{\lambda,k}[\tilde\mu]$ denotes the corresponding 
weight space. 
Since $V_{\lambda,k}$ is a highest weight module of highest
weight $(\lambda,k,0)$, we know that for all weights 
$\widetilde\mu\in P(\lambda,k)$, the difference
$(\lambda,k,0)-\widetilde\mu$ has
to be a sum of positive roots of $\widetilde L(\g)$.

For any $mu\in P$, let us set 
$$
\widetilde\mu = \left(\mu,k+h^\vee,
   \frac{\Vert\lambda+\rho\Vert^2-\Vert\mu\Vert^2}{2(k+h^\vee)}
   \right)
$$
It follows from  equation (\ref{ch}) 
that, whenever $a_\mu\neq0$ in
equation (\ref{sum}), then $(\lambda+\rho,k+h^\vee,0)-\widetilde\mu$
has to be a sum of positive roots of $\widetilde L(\g)$. 
Furthermore, we have 
$\Vert\widetilde\mu\Vert^2=\Vert(\lambda+\rho,k+h^\vee,0)\Vert^2$ and 
$\langle\widetilde\alpha,\widetilde\mu\rangle \geq 0$ for all 
$\widetilde\alpha\in\widetilde\Pi$. Together, these observations imply
that $\widetilde\mu=(\lambda+\rho,k+h^\vee,0)$ (see e.g. \cite{PS},
Lemma 14.4.7).

Putting everything together, we have proved the following theorem.

\begin{Theorem}[Kac-Weyl character formula]\label{KW}
  The character
  $\chi_{\lambda,k}$ of the integrable highest weight module $V_{\lambda,k}$
  of highest weight $\lambda$ and level $k$ at the point 
  $(h,\tau)\in\h\times\H$, where $\H$ denotes the upper half plane in $\C$,
  is given by
  $$  
    \chi_{\lambda,k}(h;\tau) = 
    \frac{\sA\Theta_{\lambda+\rho,k+h^\vee}(h,\tau)}
    {\F(h,\tau)}\,.
  $$
\end{Theorem}

%

\subsection{Characters for non-connected loop groups}
\label{Sec:twchar}
In this section  we want to derive an analogue of the Kac-Weyl
character formula in the case that $G$ not simply connected. As
always, let $G=\widetilde G/Z$ where $Z=\langle c\rangle$ is a cyclic
subgroup of the center of the simply connected Lie group 
$\widetilde G$.
Let $k$
be a multiple of the basic level of $G$, and let $\lambda\in P^k_+$ be
invariant under the action of $\sigma_c$ on  $P^k_+$. Let
$\chi_{\lambda,k}$ denote character of the $\widetilde
L^{k_b}(G)$--module $V_{\lambda,k}$. One checks directly, that under
the action of
$\sigma_c$ on $V_{\lambda,k}$, the weight space 
$V_{\lambda,k}[\widetilde\mu]$
is mapped to   $V_{\lambda,k}[\sigma_c(\widetilde\mu)]$. 
Recall the definition
of $H_0$ from the beginning of section \ref{theDE}. By definition of 
the character $\chi^{\sigma_c}_{\lambda,k}$ (see section
\ref{charsec}),
we can write
\begin{equation}\label{ch2}
  \chi^{\sigma_c}_{\lambda,k}( H_0) = 
   \underset{\sigma_c(\widetilde\mu)=\widetilde\mu}
    {\sum_{\widetilde\mu\in P(\lambda,k)}}
   Tr(\sigma_c)|_{V_{\lambda,k}[\widetilde\mu]}
   e^{2\pi i\langle\widetilde\mu,H_0\rangle}\,.
\end{equation}

On the other hand, $W_0$-anti-invariance of
$\chi^{\sigma_c}_{\lambda,k}$, and the differential equation from
Theorem \ref{thmdgl} together with Proposition \ref{thetabasis2} imply
that we can write
\begin{equation}\label{sum2}
  \chi^{\sigma_c}_{\lambda,k}( H_0)\F_{w_c}(H_0) = 
    \sum_{\mu\in \Lambda(T_{w_c})^*_{k_b}\cap (k+h^\vee)\mathfrak{a}_{w_c}}
       a_\mu f_\mu\sA_0\Theta_{\mu,k+h^\vee}
\end{equation}
where $a_\mu$ is a constant depending only on $\mu$, and
$$
  f_\mu=\exp\left(\frac{2\pi i \tau}{2(k+h^\vee)}
  (\langle\lambda+\rho_{w_c},\lambda+\rho_{w_c}\rangle - 
  \langle\mu,\mu\rangle)\right)\,.
$$

Let us set
$$
 \widetilde\mu = 
  \left(\mu,k+h^\vee,
   \frac{\Vert\lambda+\rho_{w_c}\Vert^2-\Vert\mu\Vert^2}{2(k+h^\vee)}\right)
$$
Similar to the simply connected case, 
it follows from equations (\ref{ch2}) and (\ref{sum2}) that whenever
$a_\mu\neq0$ in equation (\ref{sum2}), then 
$(\lambda+\rho_{w_c},k+h^\vee,0) - \widetilde\mu$
has to be a $\sigma_c$-invariant sum of positive roots of $\widetilde L(\g)$.
Furthermore, we have 
$$ 
  \Vert(\lambda+\rho_{w_c},k+h^\vee,0)\Vert^2=\Vert\widetilde\mu\Vert^2
$$
and $\langle\widetilde\alpha,  \widetilde\mu\rangle\geq0$ for all
$\widetilde\alpha\in\widetilde\Pi$. As in the simply connected case, these
observations imply 
$\lambda+\rho_{w_c}=\mu$. So we have proved the following theorem.

\begin{Theorem}\label{KW2}
  The character
  $\chi_{\lambda,k}$ of the integrable highest weight module $V_{\lambda,k}$
  of highest weight $\lambda$ and level $k$ at the point 
  $(\sigma_c H_0)$ is given by 
  $$  
    \chi^{\sigma_c}_{\lambda,k}(H_0) = 
    \frac{\sA_0\Theta_{\lambda+\rho_{w_c},k+h^\vee}(H_0)}
    {\F_{w_c}(H_0)}\,.
  $$
\end{Theorem}

\subsection{The root system $\widetilde\Delta_{\sigma_c}$}\label{sec:del}

Finally, let us take a closer look at the denominator $\F_{w_c}(H_0)$
appearing in Theorem \ref{KW2}. We have to distinguish two
cases. First, let us suppose that for any simple root
$\widetilde\alpha\in\widetilde\Pi$, the roots $\widetilde\alpha$ and
$\sigma_c(\widetilde\alpha)$ are not connected in the Dynkin diagram
of $\widetilde\Delta$. In this case, one can easily show that
$\sigma_c(e_{\widetilde\alpha})=e_{\sigma_c(\widetilde\alpha)}$, so that
$s(\widetilde\alpha)=1$ for all real roots
$\alpha\in\widetilde\Delta^{re}$.
For any root $\widetilde\alpha\in\widetilde\Delta$ let us denote by
$\widetilde\alpha_{\sigma_c}$ its restriction to the subspace
$\widetilde\h^{\sigma_c}\subset\widetilde\h$. Then the set
$$
  \{m_{\widetilde\alpha}\widetilde\alpha_{\sigma_c}
    ~|~\widetilde\alpha\in\widetilde\Delta^{re}\}
$$
is the set of real roots of an affine root system which we will denote
by $\widetilde\Delta_{\sigma_c}$. 

Now suppose that there exists some 
$\widetilde\alpha\in\widetilde\Delta$ such
that $\widetilde\alpha\neq\sigma_c(\widetilde\alpha)$ are not
orthogonal. Since we have excluded the case that $\sigma_c$ is the
order $n+1$ automorphism of the extended Dynkin diagram of
$\textnormal{A}_n$, we can assume that
$\sigma_c^2(\widetilde\alpha)=\widetilde\alpha$. In this case, one can
show that 
$$
 \sigma_c(e_{\widetilde\alpha}) =
   -1^{\textnormal{ht}(\widetilde\alpha)+1}e_{\sigma_c(\widetilde\alpha)}\,,
$$
where $\textnormal{ht}(\widetilde\alpha)$ denotes the height of the
root $\widetilde\alpha$ with respect to the basis $\widetilde\Pi$ (see
e.g.\cite{K}, 7.10.1 for the case of finite root systems). Looking at
the ecplicit form of the automorphism $\sigma_c$, we see that if there
exists a simple root $\widetilde\alpha$ such that
$\widetilde\alpha\neq\sigma_c(\widetilde\alpha)$ 
are not orthogonal, then we have 
$\widetilde\alpha\neq\sigma_c(\widetilde\alpha)$ for all 
simple roots of $\widetilde\Delta$.
Therefore, if
$\sigma_c(\widetilde\alpha)=\widetilde\alpha$, then there exists some
$\widetilde\beta\in\widetilde\Delta$ with
$\sigma_c(\widetilde\beta)\neq\widetilde\beta$ and
$\widetilde\beta+\sigma_c(\widetilde\beta)=\widetilde\alpha$. Hence, in
this case
$\textnormal{ht}(\widetilde\alpha)$ is necessarily even so that
$s(\widetilde\alpha)=-1$ whenever
$\sigma_c(\widetilde\alpha)=\widetilde\alpha$.
As above, the set
\begin{multline*}
  \{m_{\widetilde\alpha}\widetilde\alpha_{\sigma_c}~|~
    \widetilde\alpha\in\widetilde\Delta,
    ~\sigma_c(\widetilde\alpha)\neq\widetilde\alpha,
   ~\text{and}~\widetilde\alpha
    ~\text{and}~\sigma_c(\widetilde\alpha)~\text{are orthogonal}\}~\cup~\\
  \{2m_{\widetilde\alpha}\widetilde\alpha_{\sigma_c}~|~
    \widetilde\alpha\in\widetilde\Delta,
    ~\sigma_c(\widetilde\alpha)\neq\widetilde\alpha,
   ~\text{and}~\widetilde\alpha
    ~\text{and}~\sigma_c(\widetilde\alpha)~\text{are not orthogonal}\}
\end{multline*}
is the set of real roots of an affine root system which we also denote
by $\widetilde\Delta_{\sigma_c}$. In both cases, the smallest
positive imaginary root of
$\widetilde\Delta_{\sigma_c}$ is given by 
$\delta_{\sigma_c}=\ord(\sigma_c)\delta$,
where $\delta$ denotes the smallest positive 
imaginary root of $\widetilde\Delta$.
We will list the types of $\widetilde\Delta_{\sigma_c}$ in the end of
this paper.

Now, using $(1+x)(1-x)=(1-x^2)$, we can write 
$$
  \prod_{[\widetilde\alpha]\in(\widetilde\Delta^{re}_+)^{\sigma_c}}
   \big(1-s(\widetilde\alpha)
    e^{-2\pi i m_{\widetilde\alpha}\widetilde\alpha(H_0)}\big) = 
  \prod_{\widetilde\alpha\in\widetilde\Delta^{re}_{\sigma_c+}}
   \big(1-e^{-2\pi \widetilde\alpha(H_0)}\big)
$$
From this, we see that up to the factor
$$
  f(q) = \frac{\prod_{n=1}^\infty 
    \prod_{j=1}^{\ord(w_c)}(1-\epsilon^{-j}q^n)^{\dim(\h_j)}}
  {\prod_{n=1}^\infty
   (1-q^{\ord(\sigma_c)n})^{\text{mult}(n\delta_{\sigma_c})}}\,,
$$
the function $\F_{w_c}$ can be identified with the Kac-Weyl
denominator corresponding to the affine 
root system $\widetilde\Delta_{\sigma_c}$.
Similarly, the group $\widetilde
W_{\sigma_c}=W_0\ltimes\Lambda(T_{w_c})$ is isomorphic to the Weyl
group of the root system $\widetilde\Delta_{\sigma_c}$. So the
character $\chi_{\lambda,k}$ closely resembles a character of an
irreducible highest  weight module of an affine  Lie algebra
corresponding to the root system $\widetilde\Delta_{\sigma_c}$.
\begin{Remark}
  Since we are mainly interested in the case of loop groups based on
  connected but not necessarily simply connected Lie groups $G$, we
  restricted our attention to automorphisms $\sigma_c$ of the extended
  Dynkin diagram of $\Delta$ which are associated to elements of the
  center of the universal cover $\widetilde G$ of $G$. The arguments
  in this section can be extended to the case of 
  characters of loop groups
  based on non connected Lie groups. In this case, one has to consider
  the full automorphism group of the extended Dynkin diagram of
  $\Delta$. 
\end{Remark}
\begin{Remark}  
  If $\widetilde\Delta$ is the root system of an affine Lie algebra
  $\widetilde L(\g)$, and $\sigma$ is an automorphism of the Dynkin
  diagram of $\widetilde\Delta$ 
  then the affine Lie algebra corresponding to the
  root system $\widetilde\Delta_{\sigma}$ is often called the orbit
  Lie algebra corresponding to $\widetilde L(\g)$ and $\sigma$.
  The appearance of the root system $\widetilde\Delta_{\sigma}$ 
  was also realized in \cite{FRS} \cite{FSS}, where 
  for an outer automorphism $\sigma$ of a (generalized) Kac-Moody Lie
  algebra $\g$, the $\sigma$-twisted characters of highest weight 
  representations are calculated. It turns out that these twisted
  characters can be identified with untwisted characters of the orbit
  Lie algebra corresponding to $\g$ and $\sigma$. 
\end{Remark}


\section{Appendix: Some data on affine root systems and their automorphisms}
The following table lists some data corresponding to the
non-simply connected Lie groups. The basic levels of non-simply
connected Lie groups have been calculated in \cite{T}.
See also \cite{FSS} for a list of the root systems
$\widetilde\Delta_{\sigma}$ for general automorphisms $\sigma$
of the Dynkin diagram
of $\widetilde\Delta$. The notation for affine root systems
in the table below is the same as in \cite{K}.

\bigskip

\begin{center}
\renewcommand {\arraystretch}{1.35}
\begin{tabular}{|| l l | l | l | l | l | l ||}\hline
  $\widetilde G$ &   & $\langle c\rangle$ & 
  $G=\widetilde G/\langle c\rangle$            
 & $k_b$ &$\widetilde \Delta$&$\widetilde\Delta_{\sigma_c}$ \\\hline

 SL$_{n}$  &   $n\geq 2$     &   $\Z_r$     &     &   smallest $k$ with 
  & A$_{n-1}^{(1)}$  &   A$_{n/r-1}^{(1)}$ if $r\neq n$         \\
  & &  &  &  $\frac{n(n-1)}{r^{2}}k\in\Z$ & & $\emptyset$ if $r=n$\\ \hline
 
 Spin$_{2n+1}$  &   $n\geq 2$     &   $\Z_2$ & SO$_{2n+1}$   &  1
  & B$_{n}^{(1)}$  &   $A_{2(n-1)}^{(2)}$         \\\hline

 Sp$_{4n}$  &   $n\geq 1$     &   $\Z_2$ &    & $1$
  & C$_{2n}^{(1)}$  &   $A_{2n}^{(2)}$         \\\hline

 Sp$_{4n+2}$  &   $n\geq 1$     &   $\Z_2$ &    &  $2$
  & C$_{2n+1}^{(1)}$  &   $C_{n}^{(1)}$         \\\hline

 Spin$_{4n}$  &   $n\geq 2$ &  $\Z_2^0$ & SO$_{4n}$ & 1  
 &D$_{2n}^{(1)}$ &  C$_{2n-2}^{(1)}$     \\\hline

 Spin$_{4n}$  &   $n\geq 2$     &  $\Z_2^\pm$ &  &$1$ if $n$ even   
 & D$_{2n}^{(1)}$& B$_n^{(1)}$         \\

 &   &   &  & $2$ if $n$ odd 
 & &          \\\hline

 Spin$_{4n+2}$  &   $n\geq 2$     & $\Z_2$  & SO$_{4n+2}$ &  1
 &  D$_{2n+1}^{(1)}$ &  C$_{2n-1}^{(1)}$        \\\hline

 Spin$_{4n+2}$  &   $n\geq 2$     & $\Z_4$  & PSO$_{4n+2}$ &  4
 & D$_{2n+1}^{(1)}$ & C$_n^{(1)}$         \\\hline

 E$_{6}$  &        & $\Z_3$   &  & 3 & E$_6^{(1)}$ & G$_2^{(1)}$    \\\hline

 E$_{7}$  &        & $\Z_2$  &  &  2 & E$_7^{(1)}$ & F$_4^{(1)}$     \\\hline
\end{tabular}
\end{center}


\end{document}